\documentclass[10pt,twoside,a4paper,reqno]{amsart}
\usepackage{amscd,amsmath,amsthm,amsfonts,latexsym,amssymb}

\theoremstyle{plain}
\newtheorem{theo+}           {Theorem}      [section]
\newtheorem{prop+}           {Proposition}  [section]
\newtheorem{coro+}           {Corollary}    [section]
\newtheorem{lemm+}           {Lemma}        [section]
\newtheorem{conjecture}   {Conjecture}

\theoremstyle{definition}
\newtheorem{rema+}           {Remark}       [section]
\newtheorem{defi+}           {Definition}   [section]
\newtheorem{claim}   {Claim}
\newtheorem{problem}  {Problem}

\newenvironment{theorem}{\begin{theo+}}{\end{theo+}}
\newenvironment{proposition}{\begin{prop+}}{\end{prop+}}
\newenvironment{corollary}{\begin{coro+}}{\end{coro+}}
\newenvironment{lemma}{\begin{lemm+}}{\end{lemm+}}
\newenvironment{remark}{\begin{rema+}}{\end{rema+}}
\newenvironment{definition}{\begin{defi+}}{\end{defi+}}

\begin{document}
\numberwithin{equation}{section}

\title[Maximal and inextensible polynomials]{Maximal and inextensible 
polynomials and the geometry of the spectra of normal operators}
\author[Julius Borcea]{Julius Borcea}
\dedicatory{Dedicated to Harold S.~Shapiro on his 
75th birthday}
\keywords{Sendov's conjecture, Smale's mean value conjecture, zeros and 
critical points of complex polynomials, normal operators}
\subjclass[2000]{Primary: 30C15; Secondary: 47B15}
\address{Department of Mathematics, Stockholm University, SE-106 91 
Stockholm, Sweden}
\email{julius@math.su.se}

\begin{abstract}
We consider the set $S(n,0)$ of monic complex polynomials of degree 
$n\ge 2$ having all their zeros in the closed unit disk and vanishing at 0.
For $p\in S(n,0)$ we let $|p|_{0}$ denote the distance from the 
origin to the zero set of $p'$. We determine all 0-maximal polynomials of 
degree $n$, that is, all polynomials $p\in S(n,0)$ such that 
$|p|_{0}\ge |q|_{0}$ for any 
$q\in S(n,0)$. Using a second order variational method we then 
show that although some of these polynomials are inextensible, they are not 
necessarily locally maximal for Sendov's conjecture. This invalidates the 
recently claimed proofs of the conjectures of Sendov and Smale and shows that 
the method used in these proofs can only lead to (already known) partial 
results. In the second part of the paper we obtain a 
characterization of the critical points of a complex polynomial by means of 
multivariate majorization relations. We also propose an operator 
theoretical approach to Sendov's conjecture, which we formulate in terms of 
the spectral variation of a normal operator and its compression to the 
orthogonal complement of a trace vector. Using a theorem of Gauss-Lucas type 
for normal operators, we relate the problem of locating the critical points 
of complex polynomials to the more general problem of describing the 
relationships between the spectra of normal matrices and the spectra of their 
principal submatrices.
\end{abstract}

\maketitle

\section*{Introduction}
{\allowdisplaybreaks
Let $S_{n}$ be the set of all monic complex polynomials of degree 
$n\ge 2$ having all their zeros in the closed unit disk $\bar{D}$. If 
$p\in S_{n}$ and 
$a\in Z(p)$ then the Gauss-Lucas theorem implies that 
$(a+2\bar{D})\cap Z(p')\neq \emptyset$, where $Z(p)$ and $Z(p')$ denote the 
zero sets of $p$ and $p'$, respectively. In 1958 Sendov conjectured that 
this result may be substantially improved in the following way:
\begin{conjecture}\label{send1}
If $p\in S_{n}$ and $a\in Z(p)$ then $(a+\bar{D})\cap Z(p')\neq \emptyset$.
\end{conjecture}
Sendov's conjecture is widely regarded as one of the main challenges in the 
analytic theory of polynomials. Numerous attempts to verify this conjecture 
have led to over 80 papers, but have met with limited success. We refer to 
\cite{RS}, \cite{Se} and \cite{Sh} for surveys of the results on Sendov's 
conjecture and related questions.

The set $P_n$ of monic complex polynomials of degree $n$ may be 
viewed as a metric space by identifying it with the quotient of 
$\mathbf{C}^n$ by the action of the symmetric group on $n$ elements 
$\Sigma_n$. Indeed, let $\tau:\mathbf{C}^n\rightarrow \mathbf{C}^n/\Sigma_n$ 
denote the orbit map. Let further $p(z)=\prod_{i=1}^{n}(z-z_i)$ and 
$q(z)=\prod_{i=1}^{n}(z-\zeta_i)$ be arbitrary polynomials in $P_n$ and set
$$\Delta(p,q)=
\min_{\sigma\in \Sigma_n}\max_{1\le i\le n}|z_i-\zeta_{\sigma(i)}|.$$
Then $\Delta$ is a distance function on $P_n$ which 
induces a structure of compact metric space on the set 
$S_n=\{p\in P_n:\Delta(p,z^n)\le 1\}=\tau(\bar{D}^n)$. 
Conjecture~\ref{send1} is therefore an extremum problem 
in the closed unit ball in $P_n$ for the function $d$ given by
$$d:S_{n}\rightarrow [0,2],
\quad p\mapsto d(p)=\max_{z\in Z(p)}\min_{w\in Z(p')}|z-w|.$$
Note that $d(p)$ is the same as the so-called directed (or oriented) 
Hausdorff distance from $Z(p)$ to $Z(p')$ (cf.~\cite{Se}). Since $d$ is 
obviously a continuous function it follows by compactness that 
there exists $p\in S_{n}$ such that $d(p)=\sup_{q\in S_{n}}d(q)$. A 
polynomial with this property is called {\em extremal} for Sendov's 
conjecture. In 1972 Phelps and Rodriguez proposed the following strengthened 
form of Sendov's conjecture (cf.~\cite{PR}):
\begin{conjecture}\label{send2}
If $p\in S_{n}$ is extremal for Sendov's conjecture then 
$p(z)=z^n+e^{i\theta}$ for some $\theta \in \mathbf{R}$.
\end{conjecture}

A proof of Conjecture~\ref{send2} and thereby of Sendov's conjecture was 
recently claimed in \cite{S1}. There are currently eight different versions 
of \cite{S1}, which we shall refer to as \cite[v$k$]{S1}, $1\le k\le 8$. 
The method employed in {\em loc.~cit.~}consists in studying 
the dynamics of the zeros of the derivative of a polynomial in $S_{n}$ under 
certain perturbations of the zeros of the polynomial itself. Using these 
perturbations, a notion of extensible polynomial is defined and what is 
actually claimed in \cite[v1-v2]{S1} is that a polynomial in $S_{n}$ is 
extensible unless it vanishes only on the unit circle. If true, such a result 
would imply that the polynomials in Conjecture~\ref{send2} are not only all 
the extremal polynomials for Sendov's conjecture but also that they are in 
fact all the local maxima for the function $d$. The arguments used in 
\cite[v1-v2]{S1} were subsequently modified or
replaced by completely new ones in \cite[v3-v8]{S1}. Keeping 
track of so many changes and versions is both time-consuming and 
technically challenging. Unfortunately, this development has been rather 
confusing and has led some to believe that Schmieder's proof is correct or 
that his method would eventually work after some minor adjustments. The first 
main objective of this paper is to show that this is most definitely not the 
case. In section 1 we make a detailed analysis of first order variational 
methods in general and Schmieder's method in particular. We produce concrete
counterexamples to the main claims in all eight versions of \cite{S1} (see 
Theorem~\ref{thm15} and Propositions~\ref{prop17}--\ref{v8prop1}). Moreover, 
in sections 1.3 and 2.2 we show that 
Schmieder's approach -- or indeed any approach based exclusively on first 
order variational methods -- cannot be successful. We also argue that all 
eight versions of the proof of Smale's mean value conjecture claimed in 
\cite{S2} fail for similar reasons (section 1.4).

In spite of the failure of Schmieder's approach, variational methods remain a 
natural way of dealing with 
Sendov's conjecture. The relatively few properties which are known to 
hold for locally 
maximal or indeed even extremal polynomials were all deduced by using such 
methods (see \cite{B}, \cite{M1}, \cite{M2}, \cite{SSz}). In this spirit,  
Miller proposed a slightly more general extremal problem in \cite{M1} and 
\cite{M2}. Let $\beta\in \bar{D}$ and denote by $S(n,\beta)$ the 
set of all polynomials in $S_{n}$ which have at least one zero at $\beta$. For 
$\alpha\in \mathbf{C}$ and $p\in S(n,\beta)$ let 
$$|p|_{\alpha}=\min_{w\in Z(p')}|\alpha-w|$$ 
and define the $\alpha$-{\em critical circle} to be the circle with center 
$\alpha$ and radius $|p|_{\alpha}$. If $p\in S(n,\beta)$ is such that 
$|p|_{\alpha}\ge |q|_{\alpha}$ for 
any $q\in S(n,\beta)$ then $p$ is said to be {\em maximal with respect to} 
$\alpha$ in $S(n,\beta)$. For the sake of simplicity, maximal polynomials 
with respect to $\alpha$ in 
$S(n,\alpha)$, $\alpha\in \bar{D}$, will be called $\alpha$-{\em maximal} 
throughout this paper. A compactness argument similar to the one used for 
Conjecture~\ref{send1} 
shows that maximal polynomials do exist for any $\alpha \in \mathbf{C}$ and 
$\beta\in \bar{D}$ (cf.~\cite[Proposition 2.3]{M1}). In 1984 Miller made the 
following conjecture:

\begin{conjecture}\label{miller}
If $p\in S(n,\beta)$ is maximal with respect to $\alpha$ then all the zeros 
of $p'$ lie on the $\alpha$-critical circle and, given this, as many zeros of 
$p$ as possible lie on the unit circle.
\end{conjecture}

For results pertaining to Miller's conjecture we refer to \cite{B}, 
\cite{M1}, \cite{M2} and \cite{SSz}. 
The relevance of $\alpha$-maximal polynomials in this 
context is quite clear. As shown in \cite[Proposition 2.4]{M1}, if 
$p\in S_{n}$ is an extremal polynomial for Sendov's conjecture then there 
exists $\alpha\in Z(p)$ such that $p$ is $\alpha$-maximal and 
$|p|_{\alpha}=d(p)$. 
{\em \`A priori} there may exist $\alpha\in \bar{D}$ such that if $p$ is an 
$\alpha$-maximal polynomial then $|p|_{\alpha}<d(p)$. The $\alpha$-maximal 
polynomials that satisfy $|p|_{\alpha}=d(p)$ are 
particularly interesting not only because all extremal polynomials are 
necessarily of this type -- as already mentioned above -- but also because 
they may provide potential 
candidates for local maxima for the function $d$. 
It has been known for quite some time now that if 
$|\alpha|=1$ then $z^n-\alpha^n$ is the only $\alpha$-maximal polynomial 
(cf.~\cite{Ru}). Moreover, this 
polynomial was shown to be a local maximum for $d$ (\cite{M3}, 
\cite{VZ}). These are in fact all the examples of $\alpha$-maximal 
polynomials known so far, 
as no such polynomials were found explicitly for $|\alpha|<1$. In section 2.1 
we determine all 0-maximal polynomials 
(Theorem~\ref{class0max}) and study their properties. It turns out 
that all these polynomials satisfy Miller's conjecture as well as the relation
$|p|_{0}=d(p)$. We next consider a special class of 0--maximal polynomials, 
namely rotations of the polynomial $p(z)=z^n+z$. If $n\ge 4$ then 
$p$ is locally maximal for a large class of variations of its zeros 
(Proposition~\ref{propvar4}). Furthermore, $p$ is locally maximal for the 
restriction of the function $d$ to $S(n,0)$ and it is also inextensible with 
respect to 0. These properties and the symmetrical distribution of the 
zeros and critical points of $p$ seem to suggest that if $n\ge 4$ then $p$ 
and its rotations could in fact be locally maximal for Sendov's conjecture.
The discussion in section 1 shows that first order variational methods 
are not enough for deciding whether this is true or not. In section 2.2 we use 
a second order variational method to prove that -- contrary to what one might 
expect from the aforementioned properties -- the polynomial $p$ is not 
locally maximal for Sendov's conjecture. Indeed,  
Theorems~\ref{secord4} and~\ref{secord5} show that if $n=4$ or 5 then 
$p$ is a kind of inflection point for the function $d$. We conjecture that the 
same is actually true for all degrees and also that 
polynomials of the form $z^n+e^{i\theta}$, $\theta\in \mathbf{R}$, are in fact 
all the local maxima for $d$ (Conjecture~\ref{send3}). These
results complement those obtained in section 1 and show quite clearly that the 
methods used in~\cite{S1} cannot provide successful ways of dealing with 
Sendov's conjecture in its full generality. 

So far, almost all the results on Sendov's conjecture and related questions 
were obtained by analytical arguments. As pointed out in~\cite{B}, the 
fact that $d\circ \tau$ fails to be a (logarithmically) plurisubharmonic 
function in the polydisk $\bar{D}^n$ accounts for many of the difficulties in 
studying locally maximal polynomials for Sendov's conjecture. On the other 
hand, the geometrical information contained in the 
Gauss-Lucas theorem is hardly sufficient for dealing with 
Conjectures~\ref{send1} and~\ref{send2}. This is mainly because of the 
implicit nature of the relations between the zeros and critical points of 
complex polynomials. 
Describing these relations geometrically and as explicitly as 
possible would be helpful for a great many questions in 
the analytic theory of polynomials. In section 3 we propose an operator 
theoretical interpretation and approach to Conjectures~\ref{send1} 
and~\ref{send2} (Conjecture~\ref{operator}). Moreover, we show that these 
conjectures may be viewed as part of the more general problem of describing 
the relationships between the spectra of normal matrices and the spectra of 
their principal submatrices (Problem~\ref{norm}). We also give a geometrical 
characterization of the 
critical points of complex polynomials by means of multivariate majorization 
relations. These results use a combination of operator theoretical tools and 
methods of majorization theory. Such methods were the key to Pereira's 
recent solutions to the 1947 conjecture of de Bruijn and Springer and the 
1986 conjecture of Schoenberg. Similar ideas were used by Malamud in 
\cite{Ma}, where he not only proved these same two conjectures 
-- independently and almost at the same time as 
Pereira -- but he also obtained a remarkable generalization of the de 
Bruijn-Springer conjecture (see section 3.2). The methods of
\cite{Ma} and \cite{P} seem to be particularly well suited for studying 
extremal problems for which the loci of the zeros of extremal polynomials are 
(conjectured to be) lines in the complex plane. We believe that the results 
and the setting developed in section 3 should prove useful for investigating 
the ``spectral form'' of Sendov's conjecture (Conjecture~\ref{operator}) as
well as geometrical properties of the spectra of normal matrices and their 
degeneracy one principal submatrices (Problem~\ref{norm}). 

\section{First order variational methods and the conjectures \\ of Sendov and 
Smale}

Variational methods are a natural approach to both Sendov's conjecture and 
Smale's mean value conjecture. As a matter of fact, most of the results 
concerning the general cases of these conjectures were obtained by such 
methods (see, e.~g., \cite{B}, \cite{M1}, \cite{M2}, \cite{M3}, \cite{SSz}, 
\cite{Ti}, \cite{Ty}, \cite{VZ}). In this section we make a detailed analysis 
of the proofs of these two conjectures that were recently claimed by 
Schmieder in \cite{S1} and \cite{S2}. In particular, we show that all sixteen 
versions of \cite{S1} and \cite{S2} are incorrect and also that Schmieder's 
approach -- or indeed any approach based exclusively on first order 
variational methods -- cannot be successful.

\subsection{Some inextensible polynomials}

The notion of extensible polynomial introduced in \cite[v1-v2]{S1} amounts 
to a solvability condition for a certain system of inequalities. This system 
is linear in the generic case 
when both the polynomial and its derivative have only simple zeros. The main 
claim in \cite[v1-v2]{S1} is that a polynomial in $S_{n}$ is extensible 
unless it vanishes only on the 
unit circle. In this section we produce counterexamples to this claim for 
all degrees in the 
generic case. We also show that in most of the non-generic cases the system 
of inequalities 
mentioned above is not linear. Thus, as it stands in \cite[v1-v2]{S1}, 
the notion of extensible polynomial is not well defined in these non-generic 
cases.

Although the following definition of local maximality may seem obvious, it 
is actually quite different from the one used in \cite{S1}.

\begin{definition}\label{locmaxdef}
A polynomial $p\in S_n$ is called locally maximal for Sendov's conjecture if 
it is a local maximum for the function $d$, i.~e., if there exists 
$\varepsilon>0$ such that for any $q\in S_n$ satisfying 
$\Delta(q,p)<\varepsilon$ one has $d(q)\le d(p)$.
\end{definition}
For $p\in S_{n}$ we shall use the following notations:
\begin{equation}\label{poly1}
\begin{split}
&p(z)=\prod_{i=1}^{n}(z-z_{i}),\quad a=z_{1},\quad |p|_{a}=d(p),\quad p'(z)=
n\prod_{j=1}^{n-1}(z-w_{j}),\\
&|w_{j}-a|=|p|_{a}\text{ for }1\le j\le r,\quad |w_{j}-a|>|p|_{a}\text{ for }
j\ge r+1.
\end{split}
\end{equation}
Let $h_{1},\ldots,h_{n}\in \bar{D}$, $t\in [0,1[$, and set
\begin{equation}\label{poly2}
\begin{split}
&z_{i}(t)=z_{i}(t,h_{1},\ldots,h_{n})=\frac{th_{i}+z_{i}}{1+t\overline{h}_{i}
z_{i}},\quad 1\le i\le n,\\
&q(z)=q(z;t,h_{1},\ldots,h_{n})=\prod_{i=1}^{n}(z-z_{i}(t)),\quad q'(z)=
n\prod_{j=1}^{n-1}(z-w_{j}(t)).
\end{split}
\end{equation}
Note that $q\in S_{n}$, $z_{i}(0)=z_{i}$, $1\le i\le n$, and 
$w_{j}(0)=w_{j}$, $1\le j\le n-1$. Let us assume for now that
\begin{equation}\label{poly3}
\textit{$p$ and $p'$ have only simple zeros.}
\end{equation}
By the implicit function theorem the curves $w_{j}(t)$, $1\le j\le n-1$, are 
differentiable in a neighborhood of 0. Straightforward computations lead to 
the following result (cf.~\cite[v2, Lemma 1]{S1}):

\begin{proposition}\label{prop11}
For all sufficiently small $t>0$ and $1\le j\le r$ one has
$$|w_{j}(t)-z_{1}(t)|=|p|_{a}\!\left\{1+\Re\!\left[\sum_{i=1}^{n}\left(a_{i}
(w_{j})+\overline{b_{i}(w_{j})}\,\right)h_{i}\right]t+\mathcal{O}(t^2)
\right\},$$
where 
\begin{align*}
&a_{1}(w_{j})=-\frac{1}{w_{j}-a}\!\left[1+\frac{p(w_{j})}{(w_{j}-a)^{2}
p''(w_{j})}\right],\\
&a_{i}(w_{j})=-\frac{p(w_{j})}{(w_{j}-a)(w_{j}-z_{i})^{2}p''(w_{j})},\quad 2
\le i\le n,
\end{align*}
and $b_{i}(w_{j})=-z_{i}^{2}a_{i}(w_{j})$, $1\le i\le n$.
\end{proposition}

\begin{remark}\label{firstnot}
The notation $b_{i}(w_{j})$ was used in \cite[v1-v2]{S1}. The 
coefficients $b_{i}(w_{j})$ were subsequently denoted by $a_{i}^{*}(w_{j})$ 
in \cite[v3-v5]{S1}.
\end{remark}

From Proposition~\ref{prop11} we see that if there exist $h_{1},\ldots,h_{n}
\in \bar{D}$ such that
\begin{equation}\label{syst1}
\Re\!\left[\sum_{i=1}^{n}\left(a_{i}(w_{j})+\overline{b_{i}(w_{j})}\,\right)
h_{i}\right]>0,\quad 1\le j\le r,
\end{equation}
then $d(q)\ge |q|_{z_{1}(t)}=\min_{1\le j\le r}|w_{j}(t)-z_{1}(t)|>|p|_{a}=
d(p)$ for all small $t>0$ and thus $p$ cannot be a local maximum for $d$. 
Define the following $r\times n$ matrix:
\begin{equation}\label{B}
\mathbf{B}(p)=(\beta_{ij}),\text{ where }\,\beta_{ij}=a_{j}(w_{i})+
\overline{b_{j}(w_{i})},\quad 1\le i\le r,\,1\le j\le n,
\end{equation}
and set $\vec{h}=(h_{1},\ldots,h_{n})$. Then \eqref{syst1} may be written as
\begin{equation}\label{syst2}
\Re\left(\mathbf{B}(p)\vec{h}\right)>\vec{0}.
\end{equation}
From the above discussion one concludes that if system \eqref{syst2} has a 
solution $\vec{h}\in \mathbf{C}^n$ then the polynomial $p$ cannot be a local 
maximum for the function $d$. The following definition is essentially the 
same as Definition 2 in \cite[v2]{S1}.

\begin{definition}\label{def12}
The polynomial $p$ is {\em extensible} with respect to its zero 
$a$ if system \eqref{syst2} is solvable. Otherwise, the polynomial $p$ is 
said to be {\em inextensible}.
\end{definition}

In order to see what Definition~\ref{def12} means in terms of the 
coefficients of $\mathbf{B}(p)$ let us first recall 
\cite[Definition 2.19]{M1}.

\begin{definition}\label{def13}
A complex $m\times n$ matrix $\mathbf{M}=(m_{ij})$ is {\em positively 
singular} if there exist $\mu_{1},\ldots,\mu_{m}\ge 0$, not all 0, so that 
$\sum_{i=1}^{m}\mu_{i}m_{ij}=0$ for $1\le j\le n$.
\end{definition}

The next result is a generalization of what is usually called the fundamental 
theorem of linear programming and may be found in \cite[Theorem 22.2]{Ro} 
(see also \cite[Lemma 2.20]{M1}).

\begin{theorem}\label{thm14}
Let $\mathbf{M}$ be a complex $m\times n$ matrix and let $\vec{z}\in 
\mathbf{C}^n$ denote a vector of complex unknowns. The system 
$\Re\left(\mathbf{M}\vec{z}\right)>\vec{0}$ 
has no solution if and only if $\mathbf{M}$ is positively singular.
\end{theorem} 

From Definitions~\ref{def12}--\ref{def13} and Theorem~\ref{thm14} we see 
that $p$ is extensible with respect to its zero $a$ if and only if  
$\mathbf{B}(p)$ is not positively singular. The {\em raison d'\^etre} of 
Definition~\ref{def12} is now quite clear: extensible polynomials cannot be 
local maxima for the function $d$. Actually, this same necessary criterion 
for local maximality has already been obtained and used by Miller 
in~\cite{M1}. Note though that inextensible polynomials 
need not be local maxima for $d$ either, as one can see from the following 
simple example.

\begin{lemma}\label{countlm3}
If $\theta\in \mathbf{R}$ then the polynomial $p(z)=z^3+e^{i\theta}z$ is 
inextensible with respect to 0 but it is not a local maximum for $d$.
\end{lemma}

\begin{proof}
The fact that the polynomial $p$ is not extensible with respect to 0 is a 
special case of Theorem~\ref{thm15} below. Note that $d(p)=|p|_{0}=
\frac{1}{\sqrt{3}}$ and that by a rotation we may 
assume that $p(z)=z^3-z$. Let $t\in [0,1]$ and set
$$q_{t}(z)=(z-it)(z^2-1).$$
Then $q_{t}\in S_{3}$ and elementary computations show that for all small 
$t>0$ one has
$$d(q_{t})=|q_{t}|_{it}=\sqrt{\frac{1+t^2}{3}}>d(p),$$
which proves the lemma.
\end{proof}

\begin{remark}\label{otherex}
In section 2.2 we shall construct explicit examples of inextensible 
polynomials of degree greater than three which are not local maxima for $d$.
\end{remark}

The main claim in \cite[v1-v2]{S1} is as 
follows:

\begin{claim}\label{clm1}
If $n\ge 4$ then the polynomial $p$ is extensible with respect to its zero 
$a$ unless $p$ vanishes only on the unit circle.
\end{claim}

If true, Claim~\ref{clm1} would imply that the polynomials in 
Conjecture~\ref{send2} 
are not only all the extremal polynomials for Sendov's conjecture but also 
that they are in fact all 
the local maxima for the function $d$. The following theorem shows that 
Claim~\ref{clm1} is actually false.

\begin{theorem}\label{thm15}
Let $\theta \in \mathbf{R}$,  $n\ge 3$, and $p(z)=z^n+e^{i\theta}z$. Then 
$|p|_{0}=d(p)$ and $p$ is inextensible with respect to 0.
\end{theorem}

\begin{proof}
It is enough to prove the statement for the polynomial $p(z)=z^n-z$ since 
$z^n+e^{i\theta}z=e^{-in\alpha}p(e^{i\alpha}z)$, where 
$\alpha=\frac{\pi-\theta}{n-1}$. Recall the notations of 
\eqref{poly1} and set
$$(a=)\,z_{1}=0,\,z_{j}=e^{\frac{2\pi i(j-2)}{n-1}},\,2\le j\le n,\quad
w_{k}=\left(\!\frac{1}{n}\!\right)^{\!\frac{1}{n-1}}\!
e^{\frac{2\pi i(k-1)}{n-1}},\,1\le k\le n-1.$$
Note that $|p|_{0}=d(p)$ and $r=n-1$, so that the matrix $\mathbf{B}(p)$ 
defined in \eqref{B} is actually a $(n-1)\times n$ matrix. 
Proposition~\ref{prop11} and elementary computations yield
$$a_{1}(w_{i})=-\frac{n+1}{nw_{i}} \text{ and } a_{j}(w_{i})=
\frac{w_{i}}{n(w_{i}-z_{j})^2},\quad 1\le i\le n-1,\,2\le j\le n.$$
Using the fact that for $z\in \mathbf{C}\setminus Z(p')$ one has the 
well-known identities
\begin{eqnarray*}
&&\frac{p''(z)}{p'(z)}=\sum_{i=1}^{n-1}\frac{1}{z-w_{i}}\,\text{ and}\\
&&\frac{p^{(3)}(z)p'(z)-(p''(z))^2}{(p'(z))^2}=\frac{d}{dz}\!
\left[\frac{p''(z)}{p'(z)}\right]=-\sum_{i=1}^{n-1}\frac{1}{(z-w_{i})^2}
\end{eqnarray*}
one can show that if $2\le j\le n$ then
\begin{eqnarray*}
\sum_{i=1}^{n-1}\frac{w_{i}}{(w_{i}-z_{j})^2}&=&\sum_{i=1}^{n-1}
\frac{1}{w_{i}-z_{j}}+z_{j}\sum_{i=1}^{n-1}\frac{1}{(w_{i}-z_{j})^2}\\
&=&-\frac{p''(z_{j})}{p'(z_{j})}+z_{j}\!\left[\frac{(p''(z_{j}))^2-p^{(3)}
(z_{j})p'(z_{j})}{(p'(z_{j}))^2}\right]=\frac{n}{z_{j}}.
\end{eqnarray*}
Since $\beta_{ij}=a_{j}(w_{i})+\overline{b_{j}(w_{i})}=a_{j}(w_{i})-
\overline{z_{j}^{2}a_{j}(w_{i})}$ it follows that
$$\sum_{i=1}^{n-1}\beta_{i1}=-\frac{n+1}{n}\sum_{i=1}^{n-1}\frac{1}{w_{i}}=0 
\,\text{ and } \sum_{i=1}^{n-1}\beta_{ij}=\frac{1}{z_{j}}-\bar{z}_{j}=0,
\quad 2\le j\le n,$$
which shows that $\mathbf{B}(p)$ is positively singular. By 
Theorem~\ref{thm14} the polynomial $p(z)=z^n-z$ cannot be extensible with 
respect to its zero $z_{1}=0$.
\end{proof}

In order to explain why the proof of Claim~\ref{clm1} given in 
\cite[v1-v2]{S1} fails let us recall \eqref{poly1} and the notations of 
Proposition~\ref{prop11} and define the following $r\times n$ matrix:
$$\mathbf{A}(p)=(\alpha_{ij}),\text{ where }\,\alpha_{ij}=a_{j}(w_{i}),
\quad 1\le i\le r,\,1\le j\le n.$$
In \cite[v2, Lemma 4]{S1} it is shown that if the polynomial $p$ does not 
vanish on the unit circle then $\mathbf{B}(p)$ is positively 
singular if and only if $\mathbf{A}(p)$ is positively singular. We 
reformulate this result as 
follows:

\begin{lemma}\label{lem16}
Assume that the polynomial $p$ is as in \eqref{poly1} and that it does not 
vanish on the unit circle. Then $p$ is inextensible with respect to its zero 
$a$ if and only if $\mathbf{A}(p)$ is positively singular.
\end{lemma}

The proof of Claim~\ref{clm1} given in section 6 of~\cite[v2]{S1} is 
based on the following claim, which is a 
much stronger version of Lemma~\ref{lem16}.

\begin{claim}\label{clm2}
If $n\ge 4$ and the polynomial $p$ is inextensible with respect to its zero 
$a$ then $\mathbf{A}(p)$ is positively singular.
\end{claim}

Claim~\ref{clm2} is a crucial step in section 6 of~\cite[v2]{S1} since all the 
arguments used in {\em loc.~cit.~}rely heavily on various properties of 
$\mathbf{A}(p)$. As we shall now explain, Claim~\ref{clm2} is 
false. Indeed, the polynomials in Theorem~\ref{thm15} show that 
Lemma~\ref{lem16} cannot hold without the assumption that $p$ 
does not vanish on the unit circle:

\begin{proposition}\label{prop17}
Let $\theta \in \mathbf{R}$, $n\ge 3$, and $p(z)=z^n+e^{i\theta}z$. Then the 
polynomial $p$ is inextensible with respect to 0 and 
$\mathbf{A}(p)$ is not positively singular.
\end{proposition}

The first part of the statement in Proposition~\ref{prop17} was proved in 
Theorem~\ref{thm15}. As we shall see below, the second part of this statement 
is a consequence of the following more general result:

\begin{theorem}\label{thm18}
Let the polynomial $p$ be as in \eqref{poly1}. Assume that $p$ satisfies 
\eqref{poly3} and define the following $(n-1)\times (n-1)$ matrix:
$$\mathbf{C}(p)=(\gamma_{ij}),\text{ where }\,\gamma_{ij}=(w_{i}-z_{j})^{-2},
\quad 1\le i\le n-1,\,2\le j\le n.$$
Then $\det(\mathbf{C}(p))\neq 0$.
\end{theorem}

For the proof of Theorem~\ref{thm18} we need \cite[Lemma 2.1]{B}, which we 
restate as follows:

\begin{lemma}\label{lem19}
If the polynomial $p$ is as in \eqref{poly1} and has only simple zeros 
then there exist neighborhoods $U, V\subset \mathbf{C}^n$ of the points 
$u=(a,w_{1},\ldots,w_{n-1})$ and 
$(a,z_{2},\ldots,z_{n})$, respectively, such that
$$U\ni (\alpha,\omega_{1},\ldots,\omega_{n-1})\mapsto (\zeta_{1},\zeta_{2},
\ldots,\zeta_{n})\in V$$
is an analytic function, where $\zeta_{2},\ldots,\zeta_{n}$ are the (simple) 
zeros different from $\alpha$ of the polynomial 
$n\!\int_{\alpha}^{z}\prod_{j=1}^{n-1}(w-\omega_{j})dw$ and
$\zeta_{1}=\zeta_{1}(\alpha,\omega_{1},\ldots,\omega_{n-1})\equiv\alpha$.
\end{lemma}

\begin{remark}\label{invfcn}
We note for later purposes that if the polynomial $p$ is as in \eqref{poly1} 
and satisfies~\eqref{poly3} 
then it follows from Lemma~\ref{lem19} and the inverse function theorem that 
the functions $\omega_{1},\ldots,\omega_{n-1}$ are locally analytic in 
$\zeta_{1},\zeta_{2},\ldots,\zeta_{n}$ and one has
$$\frac{\partial \omega_j}{\partial \zeta_i}\bigg|_{v}=
-\frac{p(w_j)}{(w_j-z_i)^2p''(w_j)},\quad 1\le i\le n,\,1\le j\le n-1,$$
where $v=(z_1,z_2,\ldots,z_n)$. Using these identities one can easily compute 
the coefficients $a_i(w_j)$ in Proposition~\ref{prop11}.
\end{remark}

\noindent
{\em Proof of Theorem~\ref{thm18}.} Note first that since $p$ satisfies 
\eqref{poly3} we may choose the neighborhoods $U$ and $V$ in 
Lemma~\ref{lem19} such that if 
$(\alpha,\omega_{1},\ldots,\omega_{n-1})\in U$ then the points 
$\omega_{1},\ldots,\omega_{n-1}$ are distinct and the map 
$U\ni (\alpha,\omega_{1},\ldots,\omega_{n-1})\mapsto (\zeta_{1},\zeta_{2},
\ldots,\zeta_{n})\in V$ is onto. Since $\omega_{i}$, $1\le i\le n-1$, are 
zeros of the 
logarithmic derivative of the polynomial $\prod_{j=1}^{n}(z-\zeta_{j})=
n\!\int_{\alpha}^{z}\prod_{j=1}^{n-1}(w-\omega_{j})dw$ we get
$$\sum_{j=1}^{n}\frac{1}{\omega_{i}-\zeta_{j}}=0,\quad 1\le i\le n-1.$$
By partial differentiation we obtain
\begin{equation}\label{thm9}
\begin{split}
&\sum_{j=2}^{n}\frac{\dfrac{\partial \zeta_{j}}{\partial \omega_{i}}
\bigg|_{u}}{(w_{i}-z_{j})^2}=\sum_{j=1}^{n}\frac{1}{(w_{i}-z_{j})^2}=
-\frac{p''(w_{i})}{p(w_{i})}\neq 0,\quad 1\le i\le n-1,\\
&\sum_{j=2}^{n}\frac{\dfrac{\partial \zeta_{j}}{\partial \omega_{k}}
\bigg|_{u}}{(w_{i}-z_{j})^2}=0,\quad 1\le i\neq k\le n-1,
\end{split}
\end{equation}
where $u=(a,w_{1},w_{2},\ldots,w_{n-1})$. Moreover, by \cite[Lemma 2.3]{B}
 one has
$$\frac{\partial \zeta_{j}}{\partial \omega_{k}}\bigg|_{u}=
\frac{1}{p'(z_{j})}\int_{a}^{z_{j}}\frac{p'(w)}{w-w_{k}}dw,\quad 2\le j\le n, 
1\le k\le n-1.$$
Thus, if we set
$$\delta_{jk}=-\frac{p(w_{k})}{p'(z_{j})p''(w_{k})}\int_{a}^{z_{j}}
\frac{p'(w)}{w-w_{k}}dw,\quad 2\le j\le n, 1\le k\le n-1,$$
and define the $(n-1)\times (n-1)$ matrix $\mathbf{D}(p)=(\delta_{jk})$ then 
\eqref{thm9} may be rewritten as $\mathbf{C}(p)\mathbf{D}(p)=
\mathbf{I}_{n-1}$, where $\mathbf{I}_{n-1}$ is 
the $(n-1)\times (n-1)$ identity matrix. This implies that 
$\det(\mathbf{C}(p))\neq 0$, which proves the theorem.\hfill $\Box$

\medskip

\noindent
{\em Proof of Proposition~\ref{prop17}.} As in the proof of 
Theorem~\ref{thm15}, we may assume without loss of generality that 
$p(z)=z^n-z$. Note that for this polynomial the entries 
$\alpha_{ij}=a_{j}(w_{i})$ of the 
$(n-1)\times n$ matrix $\mathbf{A}(p)$ were already computed in the proof of 
Theorem~\ref{thm15}. Let $\mathbf{A}'(p)$ denote the $(n-1)\times (n-1)$ 
matrix which is obtained from 
$\mathbf{A}(p)$ by deleting the first column. If $\mathbf{A}(p)$ were 
positively singular then the same would have to be true for $\mathbf{A}'(p)$. 
But this is impossible since $p$ 
satisfies \eqref{poly3} and by Theorem~\ref{thm18} one has 
$$\det(\mathbf{A}'(p))=\left(-\frac{1}{n}\right)^{n}
\det(\mathbf{C}(p))\neq 0.$$
Therefore $\mathbf{A}(p)$ cannot be positively 
singular.\hfill $\Box$ 

\subsection{Multiple critical points on the critical circle}

Having explained some of the weaknesses of the arguments in \cite[v1-v2]{S1} 
in the generic case, let us now examine how these arguments are affected if 
assumption \eqref{poly3} is 
removed. We use the same notations as in \eqref{poly1}--\eqref{poly2}. 
For the sake of simplicity, we assume that 
\begin{equation}\label{nongen}
\text{$p$ {\em has only simple zeros}, $r\ge 2$, {\em and} }\,w_{1}=w_{2}=
\ldots=w_{r}.
\end{equation}
Thus, the only zero of $p'$ that lies on the $a$-critical circle of $p$ is 
$w_{1}$ with multiplicity $r\ge 2$. Let $(h_{1},\ldots,h_{n})\in \bar{D}^n$ 
and set
\begin{equation}\label{nongen1}
c=\frac{p^{(r+1)}(w_{1})}{r!},\quad d=d(h_{1},\ldots,h_{n})=-p(w_{1})
\sum_{i=1}^{n}\frac{h_{i}-\overline{h}_{i}z_{i}^2}{(w_{1}-z_{i})^2}.
\end{equation}
Note that $c\neq 0$ and that there exists a closed thin set $\Omega\subset 
\bar{D}^n$ such that $d\neq 0$ if $(h_{1},\ldots,h_{n})\in 
\bar{D}^n\setminus \Omega$. Therefore, the quantities
\begin{equation}\label{nongen2}
L_{k}=L_{k}(h_{1},\ldots,h_{n}):=\exp\!\left[i\left(\!\arg\!\left(\!
\frac{p(w_{1})}{r!c}\!\right)-\frac{r-1}{r}\arg\!\left(\frac{d}{c}\right)+
\frac{2\pi k}{r}\!\right)\right]
\end{equation}
are well defined for $1\le k\le r$ and $(h_{1},\ldots,h_{n})\in \bar{D}^n
\setminus \Omega$. In sections 3-5 of \cite[v1-v2]{S1} it is shown that if 
$(h_{1},\ldots,h_{n})\in \bar{D}^n\setminus \Omega$ then there exist 
$\varepsilon>0$ and distinct curves $w_{1j}(t)$, $1\le j\le r$, 
$t\in \,\,]0,\varepsilon]$, such that
$$\frac{\partial q(z,t,h_{1},\ldots,h_{n})}{\partial z}\!\bigg|_{w_{1j}(t)}=
0\text{ for }t\in \,\,]0,\varepsilon]\text{ and }\!
\lim_{t\rightarrow 0}w_{1j}(t)=w_{1},
\quad 1\le j\le r.$$
As in Proposition~\ref{prop11} and \eqref{B}, the coefficients of the linear 
terms in the first-order Taylor expansions of $|w_{1j}(t)-z_{1}(t)|$, 
$1\le j\le r$, are viewed as the 
entries of an 
$r\times n$ matrix $\mathbf{B}(p)=(\beta_{ij})$, $1\le i\le r$, 
$1\le j\le n$. However, these entries will now depend on the  parameters 
$h_{1},\ldots,h_{n}$. Indeed, the computations in \cite[v1-v2]{S1} 
(cf., e.~g., formula (11) in \cite[v1-v2]{S1}) show that 

\begin{equation}\label{nongen3}
\begin{split}
&\beta_{ij}=a_{j}(w_{i})+\overline{b_{j}(w_{i})},\quad 1\le i\le r,
\,1\le j\le n,\text{ where}\\
&a_{j}(w_{i})=-\frac{L_{i}}{(w_{i}-a)(w_{i}-z_{j})^2},\quad
2\le j\le n,\\
&a_{1}(w_{i})=-\frac{1}{w_{i}-a}\!\left[1+\frac{L_{i}}{(w_{i}-a)^2}\right],\,
b_{j}(w_{i})=-z_{j}^2a_{j}(w_{i}).
\end{split}
\end{equation}

\begin{remark}\label{errcorr}
The error that initially appeared in 
formula (11) of \cite[v1-v2]{S1} was later corrected in formula (20) of
\cite[v3-v4]{S1} and formula (15) of \cite[v5]{S1}.
\end{remark}

From \eqref{nongen1}--\eqref{nongen3} one can see that in 
this case system \eqref{syst2} is not linear, so that 
Theorem~\ref{thm14} can no longer be 
used. This invalidates the proof given in \cite[v1-v2, section 6]{S1} for 
non-generic cases since the arguments which are used in {\em loc.~cit.~}are 
based on the assumption that system \eqref{syst2} is 
linear and they rely heavily on Theorem~\ref{thm14}. For the same reasons, 
the word ``linear'' should be removed from Definition 2 in 
\cite[v1-v2, section 5]{S1} if one would still like to have a notion of 
extensible polynomial which is at least properly defined in non-generic cases.

\subsection{First order approximations of the critical 
points: further drawbacks} 

The arguments used in \cite[v1-v2]{S1} were subsequently modified or
replaced by completely new ones in \cite[v3-v8]{S1}. Keeping 
track of so many changes and versions is both time-consuming and 
technically challenging. Unfortunately, this development has been rather 
confusing and has led some to believe that Schmieder's proof is correct or 
that his method would eventually work after some minor adjustments. This is 
most definitely not the case, as we shall now explain.

First of all, as we already saw in the previous sections, the notion of 
inextensible 
polynomial which was implicitly defined in \cite[v1-v2]{S1} is much weaker 
than the notion of locally maximal polynomial for Sendov's conjecture (see 
Definition~\ref{locmaxdef}). In 
\cite[v3-v7, Definition 2]{S1} Schmieder tries to make these two notions 
synonymous by using a new definition of local maximality.
However, this new definition is still quite different from the actual
definition 
of local maximality which was given in Definition~\ref{locmaxdef}. As a 
matter of fact, it is ambiguously formulated and inaccurate on several points.
For instance, the words ``decreasing'' and ``increasing'' in this definition 
should be interchanged. Moreover, as it stands in 
\cite[v3-v7]{S1}, Definition~2 is valid only in the generic case when the 
polynomial $p$ has 
no multiple critical points on the critical circle. Finally, this definition 
is inadequate in that it does not allow arbitrary variations of the zeros of 
$p$. Indeed, it only allows
variations of the type described in~\eqref{poly2}. These variations are too 
restrictive since boundary zeros are always sent to boundary zeros (as is 
well known, locally maximal polynomials must have at least two zeros on the 
unit circle (cf., e.~g., \cite{PR})). 

Let us now examine the arguments used in \cite[v3-v8]{S1}.

\subsubsection*{Versions 3-5 of~\cite{S1}} One of the main claims 
in~\cite[v3-v5]{S1} is 
Theorem 1, which asserts that if a polynomial $p$ has a multiple critical point
$\zeta$ on the critical circle and if $p(\zeta)\neq 0$ then $p$ cannot be 
locally maximal unless all its zeros lie on the unit circle. However, the 
proof of Theorem 1
given in~\cite[v3-v5]{S1} is not valid. Indeed, a crucial part of this proof 
is contained in the paragraph starting with 
``For such $\mathbf{h}$ we investigate the sign of the expressions 
$F_{1l}(t)$ ...'' that immediately precedes Theorem 1. This paragraph 
contains several erroneous arguments: for instance, the conditions 
``$F_{1l}(t)\le 0$ for 
all $-\varepsilon<t<\varepsilon$ and all such $l$'' are clearly wrong and 
should be replaced by  ``$\min_{1\le l\le \sigma_k}F_{1l}(t)\le 0$ for 
all $-\varepsilon<t<\varepsilon$''. This 
invalidates the concluding argument which basically says that if $p$ is 
locally maximal then 
$0=a_m(\zeta_{1l})+a_{m}^{*}(\zeta_{1l})=a_m(\zeta_{1l})
-\overline{z_{m}^2a_m(\zeta_{1l})}$ for 
$1\le l\le \sigma_k$ and so $|z_m|=1$.

Another major objection to the arguments used in \cite[v3-v5]{S1} is that 
these are based on a criterion for local maximality 
(\cite[v3-v4, Lemma 2]{S1}) which is incorrect. Indeed, this criterion 
essentially claims that a polynomial $p$ with simple critical points on the 
critical circle is locally maximal if and only if system~\eqref{syst2} is not 
solvable or, equivalently, if and only if the matrix $\mathbf{B}(p)$ defined 
in~\eqref{B} is 
positively singular (cf.~Theorem~\ref{thm14}). This amounts to saying that 
$p$ is locally maximal if and only if it is
inextensible in the sense of Definition~\ref{def12}. The example given in
Lemma~\ref{countlm3} shows 
that this is definitely wrong. As one can see from Proposition~\ref{prop11}, 
if $t>0$ is sufficiently 
small and the polynomial $q$ is as in~\eqref{poly2} then the best one can 
say is that if $p$ is inextensible then $d(q)\le d(p)+\mathcal{O}(t^2)$. 
If anything, this shows that if $p$ is inextensible then first order 
approximations of the critical points of $p$ -- like those used in all 
versions of~\cite{S1} -- are not enough for deciding whether $p$ is locally 
maximal or not. To do this one has to use higher order 
approximations of the critical points (see section 2.2 below). Thus, the 
correct version of Lemma 2 in \cite[v3-v4]{S1} should be stated as the 
following necessary criterion for local maximality: if a polynomial $p$ with 
simple critical points on the critical circle is locally maximal then 
its associated matrix $\mathbf{B}(p)$ is positively singular. As mentioned
in section 1.1, this criterion has already been used in~\cite{M1}.
Note also that a somewhat more flexible necessary criterion for 
local maximality was obtained in~\cite{B} (the latter criterion 
imposes no conditions on the critical points and assumes only
that $p$ has simple zeros on the unit circle).

\subsubsection*{Versions 6-8 of~\cite{S1}} The main 
claim of \cite[v6-v8]{S1} is Theorem 1, which may be formulated as follows:

\begin{claim}\label{v8cl1}
Let $p\in S_n$ and $z_1\in D$ be such that $p(z_1)=0$. Then there exist a 
complex number $w_0$ with $|w_0|=1$ and $p^{*}\in S_n$ such that 
$p^{*}(w_0)=0$ and $\min_{\zeta\in Z(p^{*'})}|w_0-\zeta|\ge 
\min_{\omega\in Z(p')}|z_1-\omega|$.
\end{claim}

The proof of Claim~\ref{v8cl1} given in~\cite[v6-v7]{S1} is quite different 
from the one given in~\cite[v8]{S1}. In~\cite[v6-v7]{S1} Schmieder considers 
an arbitrary polynomial $p(z)=(z-z_1)r(z)\in S_n$ such that 
$d(p)=|p|_{z_1}=\min_{\omega\in Z(p')}|z_1-\omega|$ and he defines the 
following perturbations of $p(z)$:
\begin{equation*}
\begin{split}
&Q(z)=Q(z;t,h)=\left(z-z_1(t,h)\right)r(z),\text{ where}\\
&t\in [0,1],\,h\in \partial D
\text{ and } z_1(t,h)=\frac{z_1+th}{1+t\bar{h}z_1}.
\end{split}
\end{equation*}
A careful examination shows that if true, the local arguments used in section 
3 of \cite[v6]{S1} and section 2 of \cite[v7]{S1} would actually imply that 
the following claim must be valid as well.

\begin{claim}\label{v67cl}
There exists $h\in \partial D$ such that for all sufficiently small $t>0$ one 
has 
$\min_{\zeta\in Z(Q')}|z_1(t,h)-\zeta|>\min_{\omega\in Z(p')}|z_1-\omega|$, 
that is, $d(Q)>d(p)$.
\end{claim}

We shall now construct counterexamples to Claim~\ref{v67cl} for each degree 
$n\ge 4$. To do this, recall the notations used in the proof of 
Theorem~\ref{thm15} and let
\begin{equation}\label{c-ex1}
\begin{split}
&p(z)=z^n-z=\prod_{j=1}^{n}(z-z_j)\text{ and }
p'(z)=n\prod_{k=1}^{n-1}(z-w_k),\text{ where }z_{1}=0,\\
&z_{j}=e^{\frac{2\pi i(j-2)}{n-1}},\,2\le j\le n, 
\text{ and } w_{k}=\left(\!\frac{1}{n}\!\right)^{\!\frac{1}{n-1}}\!
e^{\frac{2\pi i(k-1)}{n-1}},\,1\le k\le n-1.
\end{split}
\end{equation}
Let further $\kappa$ be a fixed positive number and 
$(\varepsilon_1,\ldots,\varepsilon_n)\in \mathbf{C}^n$ be such 
that $|\varepsilon_j|\le |\varepsilon_1|^{1+\kappa}$ for $2\le j\le n$ and set
\begin{equation}\label{c-ex2}
Q(z)=Q(z;\varepsilon_1,\ldots,\varepsilon_n)
=\prod_{j=1}^{n}(z-(z_j+\varepsilon_j)).
\end{equation}
Using the notations in~\eqref{c-ex1}--\eqref{c-ex2} we can show the following 
result:

\begin{proposition}\label{propvar4}
If $n\ge 4$ then for all sufficiently small $|\varepsilon_1|$ one has
$$\min_{\zeta\in Q'}|z_1+\varepsilon_1-\zeta|\le
\min_{\omega\in Z(p')}|z_1-\omega|
-\left[\cos\!\left(\frac{\pi}{n-1}\right)\right]|\varepsilon_1|.$$
Thus, if $\varepsilon_j=0$ for $2\le j\le n$ or, more generally, if 
$|z_j+\varepsilon_j|\le 1$ for $2\le j\le n$ then for all sufficiently small 
$|\varepsilon_1|>0$ one has $Q\in S_n$ and $d(Q)<d(p)$.
\end{proposition}
\begin{proof}
Note that $|w_k|=\left(\!\frac{1}{n}\!\right)^{\!\frac{1}{n-1}}=|p|_{z_1}=
d(p)$, $1\le k\le n-1$, and denote the zeros of $Q'(z)$ by 
$\omega_k=\omega_k(\varepsilon_1,\ldots,\varepsilon_n)$, $1\le k\le n-1$. We 
assume that these are labeled so that $\omega_k(0,\ldots,0)=w_k$, 
$1\le k\le n-1$. It follows from~\eqref{c-ex1} that if $\varepsilon_1\neq 0$  
then there exists $j\in \{1,\ldots,n-1\}$ such that 
\begin{equation}\label{realpart}
|\arg \varepsilon_1-\arg w_j|\le \frac{\pi}{n-1},\text{ so that } 
\Re\left(\frac{\varepsilon_1}{w_j}\right)\ge 
\frac{|\varepsilon_1|}{d(p)}\cos\!\left(\frac{\pi}{n-1}\right).
\end{equation}
Now using Remark~\ref{invfcn} with $v=(0,z_2,\ldots,z_n)$ and 
the computations in the proof of Theorem~\ref{thm15} together with the 
assumption
that $|\varepsilon_i|\le |\varepsilon_1|^{1+\kappa}$ for $2\le i\le n$ we get 
\begin{equation}\label{expansion}
\begin{split}
\left|\omega_j(\varepsilon_1,\ldots,\varepsilon_n)-\varepsilon_1\right|
&=\left|w_j-\varepsilon_1
+\sum_{i=1}^{n}\frac{\partial \omega_j}{\partial \zeta_i}\bigg|_{v}
\varepsilon_i+\mathcal{O}(|\varepsilon_1|^2)\right|\\
&=|w_{j}|-\frac{n+1}{n}\Re\left(\frac{\varepsilon_1}{w_j}\right)+
\mathcal{O}(|\varepsilon_1|^2).
\end{split}
\end{equation}
From~\eqref{realpart}--\eqref{expansion} and the inequality $n+1>nd(p)$ it 
follows that for $n\ge 4$ and small $|\varepsilon_1|$ one has
\begin{equation*}
\begin{split}
\left|\omega_j(\varepsilon_1,\ldots,\varepsilon_n)-\varepsilon_1\right|
&\le |w_{j}|-\left[\frac{n+1}{nd(p)}\cos\!\left(\frac{\pi}{n-1}\right)\right]
|\varepsilon_1|+\mathcal{O}(|\varepsilon_1|^2)\\
&\le d(p)-\left[\cos\!\left(\frac{\pi}{n-1}\right)\right]|\varepsilon_1|,
\end{split}
\end{equation*}
which proves the proposition.
\end{proof}
Clearly, Proposition~\ref{propvar4} contradicts Claim~\ref{v67cl} and so it 
invalidates the results of \cite[v6-v7]{S1}.

In \cite[v8]{S1} Schmieder uses a 
modified variational method and replaces the local arguments of 
\cite[v6-v7]{S1} with some global arguments. Given a polynomial 
$p(z)=(z-z_1)q(z)\in S_n$ he considers the one-parameter family of polynomials
\begin{equation}\label{v8eq1}
Q(z;u)=(z-u)q(z),\text{ where } u\in \bar{D},
\end{equation}
and studies the Riemann surface $\mathcal{R}$ that consists of the critical
points of 
$Q$ when $u$ varies in $\bar{D}$ (that is, the zeros of the equation 
$\frac{\partial}{\partial z}Q(z;u)=0$ for $u\in \bar{D}$). Note that 
the compact 
manifold $\mathcal{R}$ contains at most $2(n-1)$ branch points. 
A close examination of sections 2 and 3 of \cite[v8]{S1} shows that if 
valid, the arguments used in the proof of Claim~\ref{v8cl1} given in 
{\em loc.~cit.~}would actually imply 
that the following claim -- which is a much stronger statement than 
Claim~\ref{v8cl1} -- is also true.

\begin{claim}\label{v8cl2}
Let $p(z)=(z-z_1)q(z)$ and $p^{*}(z)=(z-w_0)q(z)$, where $q\in S_{n-1}$ and 
$z_1$ and $w_0$ are complex numbers such that $|z_1|<1$ and $|w_0|=1$. Then 
$\min_{\zeta\in Z(p^{*'})}|w_0-\zeta|\ge \min_{\omega\in Z(p')}|z_1-\omega|$.
\end{claim}

The following proposition shows that Claim~\ref{v8cl2} fails in a very strong 
way.

\begin{proposition}\label{v8prop1}
If $n\ge 5$ then there exist $z_1\in D$ and $q\in S_{n-1}$ such that for any 
complex number $w_0$ with $|w_0|=1$ one has 
$\min_{\zeta\in Z(p^{*'})}|w_0-\zeta|< \min_{\omega\in Z(p')}|z_1-\omega|$, 
where $p(z)=(z-z_1)q(z)$ and $p^{*}(z)=(z-w_0)q(z)$.
\end{proposition}

\begin{proof}
Let $w_{0}\in \mathbf{C}$ be such that $|w_0|=1$ and set $z_1=0$ and 
$q(z)=z^{n-1}-1$, so that $p(z)=z^n-z$ and 
$p^{*}(z)=(z-w_0)(z^{n-1}-1)$. For any $\omega\in Z(p')$ one has
$|z_1-\omega|=\left(\frac{1}{n}\right)^{\frac{1}{n-1}}$, so that
\begin{equation}\label{v8eq2}
\min_{\omega\in Z(p')}|z_1-\omega|=\left(\frac{1}{n}\right)^{\frac{1}{n-1}}.
\end{equation}
Since $|w_0|=1$ and $Z(q)=\{e^{\frac{2k\pi i}{n-1}}\,:\,0\le k\le n-2\}$ one 
gets from either the Schur-Szeg\"o composition theorem 
(\cite[Theorem 3.4.1d]{RS}) or the Grace-Heawood theorem 
(\cite[Theorem 4.3.1]{RS}) that
\begin{equation}\label{v8eq3}
\min_{\zeta\in Z(p^{*'})}|w_0-\zeta|\le \frac{\min_{\alpha\in Z(q)}
|w_0-\alpha|}{2\sin(\pi/n)}\le \frac{\sin\left(\pi/2(n-1)\right)}{\sin(\pi/n)}.
\end{equation}
From~\eqref{v8eq2} and~\eqref{v8eq3} we see that in order to prove the 
proposition it is enough to show that
\begin{equation}\label{v8eq4}
\frac{\sin\left(\pi/2(n-1)\right)}{\sin(\pi/n)}
< \left(\frac{1}{n}\right)^{\frac{1}{n-1}}\text{ for }n\ge 5.
\end{equation}
Numerical checking shows that~\eqref{v8eq4} is true for $n=5,6,7$ or 8, and so 
we may assume that $n\ge 9$. Since 
$0<\dfrac{\pi}{2(n-1)}<\dfrac{\pi}{n}<\dfrac{\pi}{4}$ for $n\ge 5$ and 
$x\mapsto \dfrac{\sin x}{\sqrt{x}}$ is an 
increasing function on $\left(0,\dfrac{\pi}{4}\right)$ we get that
\begin{equation}\label{v8eq5}
\frac{\sin\left(\pi/2(n-1)\right)}{\sin(\pi/n)}<\sqrt{\frac{n}{2(n-1)}}
\,\text{ for }n\ge 5.
\end{equation}
The sequence $\sqrt{2\left(1-\frac{1}{n}\right)}-n^{\frac{1}{n-1}}$ is clearly 
increasing for $n\ge 2$, so that 
$$\sqrt{2\left(1-\frac{1}{n}\right)}-n^{\frac{1}{n-1}}\ge 
\frac{4}{3}-3^{1/4}>0$$ 
whenever $n\ge 9$. This implies that the right-hand side of~\eqref{v8eq5} is 
always less than $\left(\frac{1}{n}\right)^{\frac{1}{n-1}}$ if $n\ge 9$, which 
proves~\eqref{v8eq4}. 
\end{proof}

Let us finally note that the condition $|w_0|=1$ is never really used in 
sections 2 and 3 of \cite[v8]{S1}, which is quite strange. Indeed, 
the ``blowing up and pulling back'' technique used in {\em loc.~cit.~}is 
in fact a kind of projectivization method for which no assumption on $w_0$ 
other than $|w_0|\le 1$ seems to be necessary. Proposition~\ref{v8prop1} 
shows quite clearly that such arguments cannot be valid. 

\begin{remark}\label{othererr}
We have pointed out only some of the most serious (actually, irreparable) 
errors in~\cite[v1-v8]{S1}. However, several other technical or conceptual 
mistakes appear in {\em loc.~cit.} For instance, Lemma 3 in~\cite[v3]{S1} is 
contradicted by the polynomials $p_1(z)=(z-1)^n$ and $p_2(z)=(z+1)^n$
or indeed any sufficiently small perturbations of these polynomials. Also, 
the identity ``$\log(\rho(Q(\cdot,t)))=\log |\zeta_{k}(t)-z_{1}(t)|$'' -- 
which is consistently used in~\cite[v6-v8]{S1} -- is definitely wrong and 
should be replaced by 
``$\log(\rho(Q(\cdot,t)))=\min_{k}\log |\zeta_{k}(t)-z_{1}(t)|$''.
\end{remark}

\subsubsection*{Conclusion} The failure of the method proposed in \cite{S1} 
is mainly 
due to the fact that the local arguments used in {\em loc.~cit.~}are based on 
first order approximations of the critical points of a polynomial. As we 
already explained, such arguments lead only to necessary conditions that do 
not usually provide enough information for deciding whether a given 
polynomial is locally maximal for Sendov's conjecture or not. This is 
hardly surprising since it is actually a common feature of most non-trivial 
extremal problems, which often require higher order approximations. 
A real-valued $C^1$-function $f$ of one real variable 
such that $f'$ vanishes at a point where $f$ does not have a local maximum is 
arguably the simplest example that comes to mind in this context. In section 
2.2 we shall construct explicit examples of polynomials for which second 
order approximations of their critical points are in fact the only way to 
prove that they are not locally maximal for Sendov's conjecture. It is 
therefore highly unlikely that the approach proposed in \cite{S1} could 
be made into a successful method for dealing with Sendov's conjecture in its 
full generality.

\subsection{Smale's mean value conjecture} In his 1981 work on the complexity 
of algorithms for successive 
root approximation Smale conjectured that any polynomial 
$p\in \mathbf{C}[z]$ of degree 
$n\ge 2$ with $p(0)=0$ and $p'(0)\neq 0$ satisfies the following mean value 
property 
(cf.~\cite[Problem 1E]{Sm}; see also \cite{Sm1}):
$$\min\left\{\bigg|\frac{p(w)}{p'(0)w}\bigg|\,:\,p'(w)=0\right\}\le 
\frac{n-1}{n}.$$
Reviews of the results known so far on Smale's mean value conjecture and 
related questions have appeared in  \cite{RS}, \cite{Sc} and \cite{Sh}.
A proof of this conjecture was recently claimed in \cite{S2}. The variational 
method used in {\em loc.~cit.~}is the same as the one that 
was used in \cite{S1} for Sendov's conjecture. A notion of extensible 
polynomial similar to the one defined in \cite{S1} 
is introduced and it is claimed that $p$ is extensible unless all the roots 
of $p$ other than 0 have the same absolute value. Unfortunately, the proof 
given in \cite{S2} relies 
upon the same erroneous or insufficient arguments as those described above in 
the case of Sendov's conjecture. Indeed, a close examination shows that 
except for some minor changes, the arguments used in \cite[v1-v3]{S2}, 
\cite[v4-v5]{S2}, \cite[v6-v7]{S2} and \cite[v8]{S2} are 
pretty much the same as those used in \cite[v1-v2]{S1}, \cite[v3-v5]{S1}, 
\cite[v6-v7]{S1} and 
\cite[v8]{S1}, respectively. Thus, the proofs of Smale's mean value conjecture
claimed in \cite[v1-v8]{S2} are not valid. For the same reasons, one may 
safely say that there is hardly any chance that the approach proposed in 
{\em loc.~cit.~}could be made into a successful method for dealing with 
Smale's mean value conjecture in its full generality.

On a different note, we would like to mention an interesting recent preprint 
of Tyson (\cite{Ty}) where counterexamples to Tischler's strong form of 
Smale's mean value conjecture (cf.~\cite{Ti}) were constructed for each degree 
$n\ge 5$.

\section{0-maximal polynomials, second order variational methods, and local 
maxima for Sendov's conjecture}

The polynomials $z^n-\alpha^n$ with $|\alpha|=1$ are the only concrete 
examples of $\alpha$-maximal polynomials known so far. As pointed out in the 
introduction, these are also the only known examples of local maxima  
for Sendov's conjecture. In this section we
first give a complete description of 0-maximal 
polynomials of arbitrary degree. We then use a second order 
variational method to show that -- despite all appearances -- 
these polynomials are not necessarily locally maximal for Sendov's conjecture.
We conjecture that polynomials with the latter property must in fact be 
$\alpha$-maximal with $|\alpha|=1$ (Conjecture~\ref{send3} below).

\subsection{Classification of 0-maximal polynomials}

In this section we determine all 0-maximal polynomials and study some of 
their properties. We start with a few preliminary results, the first of which 
is well known from the theory of self-inversive polynomials:

\begin{lemma}
Let $p(z)=\sum_{k=0}^{n}a_{k}z^k$ be a complex polynomial of degree $n\ge 1$. 
If all the zeros of $p$ lie on the unit circle then $a_{k}\overline{a}_{0}=
a_{n}\overline{a}_{n-k}$, $0\le k\le n-1.$
\end{lemma}

From Lemma 2.1 we deduce the following result.

\begin{lemma}
Let $q$ be a complex polynomial of degree $n\ge 2$ and $\alpha\in Z(q)$. If
$R>0$ is such that $|w-\alpha|=R$ for 
any $w\in Z(q')$ then
\begin{equation*}
nR^{2(n-1)}q'(\alpha+z)=q'(\alpha)z^{n-1}\overline{q'(\alpha+R^2\overline{z}^
{\,-1})}
\end{equation*}
for all $z\in \mathbf{C}$, that is,
\begin{equation}\label{lem2}
(n-k-1)!\,nR^{2k}q^{(k+1)}(\alpha)=k!\,q'(\alpha)\,\overline{q^{(n-k)}
(\alpha)},\quad 0\le k\le n-1.
\end{equation}
\end{lemma}
\begin{proof}
By assumption, the polynomial
$$q'(\alpha+Rz)=\sum_{k=0}^{n-1}\frac{R^kq^{(k+1)}(\alpha)}{k!}z^k$$
is of degree at least one and has all its zeros on the unit circle. The 
identities in \eqref{lem2} follow by applying Lemma 2.1 to the polynomial 
$q'(\alpha+Rz)$.
\end{proof}

We shall also need the following lemma.

\begin{lemma}
If $n$ is an integer greater than two then 
\begin{equation*}
\left\{x\in [1,n-2]\,:\, n^{\frac{2x}{n-1}}(n-x)-n(x+1)=0\right\}=
\left\{\frac{n-1}{2}\right\}.
\end{equation*}
\end{lemma}
\begin{proof}
The assertion is trivially true for $n=3$ and we may therefore assume that 
$n\ge 4$. Note that
$$(-1,n)\ni x\mapsto f(x)=\left(\frac{2x}{n-1}-1\right)\!\log n +\log\!\left(
\frac{n-x}{x+1}\right)$$
is a continuously differentiable function which satisfies $f(n-1-x)+f(x)=0$, 
so that $f\left(\frac{n-1}{2}\right)=0$. It is then easily 
seen that the lemma is in fact equivalent to the following statement:
\begin{equation}\label{lem3-1}
f(x)\neq 0 \,\text{ for any }\, x\in [1,n-2]\setminus \left\{\frac{n-1}{2}
\right\}.
\end{equation}
In order to prove \eqref{lem3-1} we set
$$\Delta_{n}=\frac{n+1}{n-1}\left(\frac{n+1}{n-1}-\frac{2}{\log n}\right)$$
and notice that
\begin{equation*}
0<\Delta_{n}<\left(\frac{n+1}{n-1}\right)^2\,\text{ if }\,n\ge 2.
\end{equation*}
It follows that $-1<x_{-}<\frac{n-1}{2}<x_{+}<n$, where $x_{\pm}=\frac{n-1}{2}
\left(1\pm \sqrt{\Delta_{n}}\right)$. An elementary 
computation yields
\begin{equation}\label{lem3-2}
f'(x)=-\frac{2\log n}{n-1}\frac{(x-x_{+})(x-x_{-})}{(n-x)(x+1)},\quad 
x\in (-1,n),
\end{equation}
which shows that the function $f$ is increasing on $(x_{-},x_{+})$ and it is 
decreasing on both $(-1,x_{-})$ and 
$(x_{+},n)$. It is not difficult to prove that
\begin{align*}
&n\in \{4,5,6\}\Rightarrow 4\log n>n+1\Rightarrow \Delta_{n}>\left(
\frac{n-3}{n-1}\right)^2\Rightarrow x_{+}>n-2\\
\intertext{and also that}
&n\ge 7\Rightarrow 4\log n<n+1\Rightarrow \Delta_{n}<\left(\frac{n-3}{n-1}
\right)^2\Rightarrow x_{+}<n-2,
\end{align*}
which together with the identity $x_{+}+x_{-}=n-1$ show that
\begin{equation}\label{lem3-3}
[1,n-2]\subset (x_{-},x_{+}) \,\text{ if }\, n\in \{4,5,6\} \,\text{ and }\,
[x_{-},x_{+}]\subset (1,n-2)\,\text{ if }\,n\ge 7.
\end{equation}
Moreover, the inequality
$$2^{n-1}\left(1+\frac{1}{n-1}\right)^{n-1}\!>n^2,\quad n\ge 4,$$
implies that
\begin{equation}\label{lem3-4}
f(n-2)=-f(1)=\frac{n-3}{n-1}\log n-\log\!\left(\frac{n-1}{2}\right)>0\,
\text{ if }\,n\ge 4.
\end{equation}
From \eqref{lem3-2}--\eqref{lem3-4} we deduce that \eqref{lem3-1} must be 
true, which proves the lemma.
\end{proof}

We are now ready to prove the main result of this section:

\begin{theorem}\label{class0max}
The 0-maximal polynomials of degree $n\ge 2$ are given by:
\begin{itemize}
\item[(i)] $z^{2m}+e^{i\theta}z$ if $n=2m$ and $m\ge 1$, where 
$\theta\in \mathbf{R}$.
\item[(ii)] $z^{2m+1}+\lambda e^{i\theta}z^{m+1}+e^{i2\theta}z$ if 
$n=2m+1$ and $m\ge 1$, where $\lambda,\theta\in \mathbf{R}$ and 
$|\lambda|\le \frac{2\sqrt{2m+1}}{m+1}$.
\end{itemize}
\end{theorem}
\begin{proof}
Let $p$ be a 0-maximal polynomial. Since $p\in S(n,0)$, we may write
\begin{eqnarray*}
&&p(z)=z^n+\sum_{k=1}^{n-1}a_{k}z^k=z\prod_{i=1}^{n-1}(z-z_{i}),
\text{ where } |z_{i}|\le 
1,\,1\le i\le n-1,\\
&&p'(z)=nz^{n-1}+\sum_{k=1}^{n-1}ka_{k}z^{k-1}=n\prod_{j=1}^{n-1}(z-w_{j}).
\end{eqnarray*}
By comparing $p$ with the polynomial $q(z):=z^{n}-z\in S(n,0)$ we get
$$\min_{1\le j\le n-1}|w_{j}|=|p|_{0}\ge |q|_{0}=\left(\!\frac{1}{n}\!\right)^
{\frac{1}{n-1}},$$
which combined with the identity $\prod_{i=1}^{n-1}z_{i}=(-1)^{n-1}p'(0)=
n\prod_{j=1}^{n-1}w_{j}$ yields 
$\prod_{i=1}^{n-1}|z_{i}|=n\prod_{j=1}^{n-1}|w_{j}|\ge 1$. Since $|z_{i}|
\le 1$ for $1\le i\le n-1$ we deduce that
\begin{equation}\label{thm2-1}
|z_{i}|=1,\,1\le i\le n-1,\,\text{ and }\,|w_{j}|=\left(\!\frac{1}{n}\!
\right)^{\frac{1}{n-1}},\,1\le j\le n-1.
\end{equation}
It follows in particular that $|a_{1}|=|\prod_{i=1}^{n-1}z_{i}|=1$, hence
\begin{equation}\label{thm2-2}
a_{1}=e^{i2\theta}
\end{equation}
for some $\theta\in \mathbf{R}$. From \eqref{thm2-1} and Lemma 2.1 applied 
to the polynomial $z^{-1}p(z)$ we get
$$\overline{a}_{n-k}=e^{-i2\theta}a_{k+1},\quad 1\le k\le n-2,$$
while \eqref{thm2-1} and Lemma 2.2 applied to the polynomial $p$ with 
$\alpha=0$ and $R=\left(\!\frac{1}{n}\!\right)^{\frac{1}{n-1}}$ 
imply that
$$\overline{a}_{n-k}=\frac{n(k+1)}{n-k}\left(\!\frac{1}{n}\!\right)^{
\frac{2k}{n-1}}e^{-i2\theta}a_{k+1},\quad 1\le k\le n-2.$$
Thus
\begin{equation}\label{thm2-3}
\left[n^{\frac{2k}{n-1}}(n-k)-n(k+1)\right]a_{k+1}=0,\quad 1\le k\le n-2,
\end{equation}
and then by Lemma 2.3 we get that $a_{k}=0$ for $2\le k\le n-1$ if $n$ 
is even. 

Let now $n=2m+1,\,m\ge 1$. In this case Lemma 2.3 and \eqref{thm2-3} imply 
that $a_{k}=0$ for $2\le k\le n-1,\,k\neq m+1$, so that
\begin{equation}\label{thm2-4}
p(z)=z^{2m+1}+a_{m+1}z^{m+1}+e^{i2\theta}z,
\end{equation}
where $\theta$ is as in \eqref{thm2-2}. From \eqref{thm2-1} and 
\eqref{thm2-4} we obtain
$$\sum_{i=1}^{2m}z_{i}^m+2ma_{m+1}+e^{i2\theta}\sum_{i=1}^{2m}
\overline{z}_{i}^m=0.$$
Note that by Newton's identities one has actually that
$$\sum_{i=1}^{2m}z_{i}^m=-ma_{m+1}+\sum_{k=1}^{m-1}a_{m+k+1}\bigg(\!
\sum_{j=1}^{2m}z_{j}^{k}\!\bigg)=-ma_{m+1}.$$
We deduce from the last two formulas that $a_{m+1}=e^{i2\theta}
\overline{a}_{m+1}$, which together with \eqref{thm2-4} implies that the 
polynomial $p$ must be of the form
$$p(z)=z^{2m+1}+\lambda e^{i\theta}z^{m+1}+e^{i2\theta}z,\,\text{ where }\,
\lambda\in \mathbf{R}.$$
An elementary computation shows that if $\theta,s,t\in \mathbf{R}$ with $t>0$ 
then the roots of the equation $x^2+se^{i\theta}x+te^{i2\theta}=0$ 
cannot have the same absolute value unless $s^2\le 4t$. By \eqref{thm2-1} all 
the zeros of the polynomial
$$p'(z)=(2m+1)z^{2m}+(m+1)\lambda e^{i\theta}z^{m}+e^{i2\theta}$$
have the same absolute value. Using the substitution $x=z^m$ and the 
above-mentioned result on second degree equations we see that this 
cannot happen unless $|\lambda|\le \frac{2\sqrt{2m+1}}{m+1}$.

To summarize, we have shown that if $p\in S(n,0)$ then $|p|_{0}<\left(\!
\frac{1}{n}\!\right)^{\frac{1}{n-1}}$ unless $p$ is of the form (i) or 
(ii) according to the parity of $n$. On the other hand, it is easily checked 
that a polynomial of the form (i) or (ii) satisfies 
$|p|_{0}=\left(\!\frac{1}{n}\!\right)^{\frac{1}{n-1}}$, which completes the 
proof of the theorem.
\end{proof}

\begin{remark}\label{remmax1}
Let $p$ be a 0-maximal polynomial of degree $n\ge 2$. By 
Theorem~\ref{class0max} all the 
critical points of $p$ lie on the 0-critical circle $|z|=\left(\!
\frac{1}{n}\!\right)^{\frac{1}{n-1}}$ and all the zeros of $p$ except 0 lie 
on the unit circle. Thus, all 0-maximal polynomials satisfy Miller's 
conjecture (Conjecture~\ref{miller}). Note also that if $p$ is 0-maximal then 
it has only simple zeros. The same is true for $p'$ except 
when $n$ is odd and $p$ is of the form (ii) with 
$|\lambda|=\frac{4\sqrt{n}}{n+1}$, in which case $p'$ has $\frac{n-1}{2}$ 
distinct zeros each of multiplicity two.
\end{remark}

\begin{corollary}\label{cormax}
If $p$ is a 0-maximal polynomial of degree $n\ge 2$ then $d(p)=|p|_0$.
\end{corollary}

\begin{proof}
Let $z\in Z(p)\setminus \{0\}$, so that $|z|=1$ by Remark~\ref{remmax1}. 
Denote by $H_z$ the closed half-plane which contains $z$ and is bounded by 
the line passing through $\frac{z}{2}$ which is orthogonal to the 
segment $[0,z]$. Since $p(z)=0=p(0)$ it follows from a well-known 
consequence of the Grace-Heawood theorem that there exists 
$w\in Z(p')\cap H_z$ (see, e.~g., the supplement to 
Theorem 4.3.1 in \cite{RS}). By 
Remark~\ref{remmax1} again one has that $|w|=|p|_0=\left(\!
\frac{1}{n}\!\right)^{\frac{1}{n-1}}\ge \frac{1}{2}$. It is then geometrically 
clear that $|z-w|\le |p|_0$, so that 
$\min_{\omega\in Z(p')}|z-\omega|\le |p|_0$.
\end{proof}  

\begin{remark}\label{remmax3}
Theorem~\ref{class0max} and Corollary~\ref{cormax} do not automatically imply 
that 0-maximal polynomials are global maxima for the restriction of the 
function $d$ to the subset $S(n,0)$ of $S_n$. For this one would have to 
show that $d(p)\le \left(\!
\frac{1}{n}\!\right)^{\frac{1}{n-1}}$ whenever $p\in S(n,0)$, which we 
actually conjecture
to be true. Note that the best estimate known so far for $p\in S(n,0)$ is 
$d(p)<1$ (see~\cite[Theorem 7.3.6]{RS}).
\end{remark}

\begin{remark}\label{remmax2}
It is not difficult to see that that $\max_{p\in S(n,1)}|p|_{1}=1$ 
(cf.~\cite{Ru}). On the other hand, Theorem~\ref{class0max} shows that 
$\max_{p\in S(n,0)}|p|_{0}$ 
increases to $1$ as $n\rightarrow \infty$. It would be interesting to know 
whether there exist $\alpha\in D$ 
and $c\in (0,1)$ such that $\max_{p\in S(n,\alpha)}|p|_{\alpha}\le c$ for any 
$n\ge 1$.
\end{remark}

To end this section let us point out that Theorem~\ref{class0max} solves 
in fact the following more general extremal problem 
concerning the distribution of zeros and critical points of complex 
polynomials:

\begin{problem}\label{pb1}
Let $n\ge 2$, $a\in \mathbf{C}$, and $R>0$. Find the largest constant 
$\rho:=\rho(a,n,R)$ with the property that for any complex polynomial $p$ of 
degree $n$ satisfying $p(a)=0$ and $\min_{w\in Z(p')}|w-a|\ge R$ one has 
$\max_{z\in Z(p)}|z-a|\ge \rho$.
\end{problem}

It is easy to see that $\rho$ is invariant under translations in the complex 
plane so that it does not depend on $a$. Moreover, the Gauss-Lucas theorem 
clearly implies that $\rho >R$. 
On the other hand, by considering $p(z)=(z-a)^n+nR^{n-1}(z-a)$ we see that 
$\rho \le Rn^{\frac{1}{n-1}}$. Essentially the same computations as in the 
proof of Theorem~\ref{class0max} show that 
one has in fact $\rho=Rn^{\frac{1}{n-1}}$:

\begin{theorem}
With the notations of Problem 1 one has $\rho=Rn^{\frac{1}{n-1}}$ and this 
value is attained only for the following polynomials:
\begin{itemize}
\item[(i)] $(z-a)^{2m}+2mR^{2m-1}e^{i\theta}(z-a)$ if $n=2m$ and $m\ge 1$,
where $\theta\in \mathbf{R}$.
\item[(ii)] $(z-a)^{2m+1}+\lambda \sqrt{2m+1}R^{m}e^{i\theta}(z-a)^{m+1}+
(2m+1)R^{2m}e^{i2\theta}(z-a)$ if $n=2m+1$ and $m\ge 1$, where 
$\lambda,\theta\in \mathbf{R}$ and $|\lambda|\le \frac{2\sqrt{2m+1}}{m+1}$.
\hfill $\Box$
\end{itemize}
\end{theorem}

\subsection{Second order approximations of the critical points and local 
maxima for Sendov's conjecture}

As shown in Lemma~\ref{countlm3}, the polynomial $z^3+z$ and its rotations are 
inextensible with respect to 0 but they are not local maxima for the function 
$d$. On the other hand, all the properties of the polynomial $p(z)=z^n+z$ 
that we discussed so far seem to suggest that if $n\ge 4$ then $p$ and its 
rotations could in 
fact be locally maximal for Sendov's conjecture. Indeed, 
Theorem~\ref{class0max} combined 
with the fact that $|p|_{\zeta}<|p|_0=d(p)$ for $\zeta\in Z(p)\setminus \{0\}$ 
implies that the polynomial $p$ is locally maximal for the restriction of the 
function $d$ to $S(n,0)$. Moreover, $p$ is locally maximal for variations such 
as those described in Proposition~\ref{propvar4}. In addition to that, the 
zeros of $p$ 
and its critical points are symmetrically distributed and satisfy 
Conjecture~\ref{miller}. Finally, by Theorem~\ref{thm15} the polynomial $p$ 
is inextensible with respect to 0. We shall now use a second order variational 
method to show that -- contrary to what one might expect from the 
aforementioned properties -- the polynomial $p$ is not locally maximal for 
Sendov's conjecture. (The discussion in section 1 clearly shows that first 
order variational methods are 
not enough for deciding whether $p$ is a local maximum for $d$ or not.) 
Theorems~\ref{secord4} and~\ref{secord5} below show that if $n=4$ or 5 then 
$p$ is a kind of inflection point for the function $d$. We conjecture that the 
same is actually true for all degrees (Remark~\ref{gencase}) and also that 
polynomials of the form $z^n+e^{i\theta}$, $\theta\in \mathbf{R}$, are in fact 
all the local maxima for the function $d$ (Conjecture~\ref{send3}). These
results complement those obtained in section 1 and show quite clearly that the 
methods used in~\cite{S1} cannot provide successful ways of dealing with 
Sendov's conjecture in its full generality. 

Let $a\in [0,1]$ and set
\begin{equation}\label{var1}
\begin{split}
&r=\left|z^4+z\right|_0=\sqrt[3]{\frac{1}{4}},
\quad \alpha_1=\frac{3\sqrt{3}r}{2-3r},\quad
\alpha_2=-\frac{\sqrt{3}\left[(3r+2)^2+4\right]}{2(3r-2)^2},\\
&\zeta(a)=\exp{\left[i\!\left(\frac{\pi}{3}+\alpha_1 a
+\alpha_2 a^2\right)\right]}.
\end{split}
\end{equation}
We use the quantities in~\eqref{var1} in order to construct certain continuous 
deformations of the polynomial $z^4+z$. For sufficiently small $a>0$ these  
perturbations give rise to a one-parameter 
family of polynomials in $S_4$ with the following interesting property:

\begin{theorem}\label{secord4}
Let $\zeta(a)$ be as in~\eqref{var1} and define the polynomial 
\begin{equation}\label{var2}
p_a(z):=(z-a)(z+1)\left[z-\zeta(a)\right]\!
\left[z-\overline{\zeta(a)}\right]\in S_4.
\end{equation}
Then for all sufficiently small $a>0$ one has 
$\Delta(p_a,z^4+z)=\mathcal{O}(a)$ and 
$$d(p_a)=|p_a|_a=r+Ca^2+\mathcal{O}(a^3),\,\text{ where }\,
C=\frac{3}{4r(2-3r)}\approx 10.81154938.$$
In particular, the polynomial $z^4+z$ is not locally maximal for 
Sendov's conjecture.
\end{theorem}

\begin{proof}
Let us first describe the steps we took in order to arrive at the 
quantities in~\eqref{var1} and the one-parameter family of polynomials 
defined in~\eqref{var2}. We start by introducing six real parameters 
which we denote by $x_i$, $y_i$, $1\le i\le 3$, and we define the auxiliary 
polynomials
\begin{eqnarray*}
P_a(z;x_1,x_2,x_3)\!&=&\!(z-a)\left[z-z_1(a;x_3)\right]
\left[z-z_2(a;x_1,x_2)\right]
\!\left[z-\overline{z_2(a;x_1,x_2)}\right]\\
&=&\!z^4+\sum_{j=0}^{3}b_j(a;x_1,x_2,x_3)z^j,\\
Q_a(z;y_1,y_2,y_3)\!&=&\!4\!\left[z-\omega_1(a;y_3)\right]
\left[z-\omega_2(a;y_1,y_2,y_3)\right]
\!\left[z-\overline{\omega_2(a;y_1,y_2,y_3)}\right]\\
&=&\!4z^3+\sum_{k=0}^{2}c_k(a;y_1,y_2,y_3)z^k,
\end{eqnarray*}
where
\begin{equation}\label{omegas}
\begin{split}
&z_1(a;x_3)=-1+x_3 a^2,\quad z_2(a;x_1,x_2)=\exp{\left[i\!\left(\frac{\pi}{3}
+x_1 a+x_2 a^2\right)\right]},\\
&\omega_1(a;y_3)=a-\left(r+y_3 a^2\right),\\
&\omega_2(a;y_1,y_2,y_3)=
a+\left(r+y_3 a^2\right)
\exp{\left[i\!\left(\frac{\pi}{3}+y_1 a+y_2 a^2\right)\right]}.
\end{split}
\end{equation}
The idea is now to investigate whether the parameters 
$x_i$ and $y_i$, $1\le i\le 3$, may be chosen so that the following 
conditions are satisfied:
\begin{itemize}
\item[(i)] $x_3\ge 0$ and $y_3>0$;
\item[(ii)] for all sufficiently small $a>0$ one has
$$\max_{|z|\le 1}\big|P'_a(z;x_1,x_2,x_3)-Q_a(z;y_1,y_2,y_3)\big|
=\mathcal{O}(a^3),$$
\end{itemize}
where $P'_a(z;x_1,x_2,x_3)$ denotes the derivative of $P_a(z;x_1,x_2,x_3)$ 
with respect to $z$. To do this, we expand the coefficients of 
$P'_a(z;x_1,x_2,x_3)$ and $Q_a(z;y_1,y_2,y_3)$ into their MacLaurin series 
so as to get second order approximations of these coefficients (with an error 
of $\mathcal{O}(a^3)$). Let $\tilde{b}_j(a;x_1,x_2,x_3)$, $1\le j\le 3$, and 
$\tilde{c}_k(a;y_1,y_2,y_3)$, $0\le k\le 2$, be the resulting second degree 
MacLaurin polynomials in the variable $a$ for the coefficients 
$b_j(a;x_1,x_2,x_3)$, $1\le j\le 3$, and $c_k(a;y_1,y_2,y_3)$, $0\le k\le 2$, 
respectively. Then we may write
\begin{multline}
(m+1)\tilde{b}_{m+1}(a;x_1,x_2,x_3)-\tilde{c}_{m}(a;y_1,y_2,y_3)\\
=\sum_{n=0}^{2}d_{mn}(x_1,x_2,x_3,y_1,y_2,y_3)a^n,\quad 0\le m\le 2,
\end{multline}
where $d_{mn}(x_1,x_2,x_3,y_1,y_2,y_3)$, $0\le m,n\le 2$, are 
real polynomials in the variables $x_1,x_2,x_3,y_1,y_2,y_3$.
Clearly, any solution $(x_1,x_2,x_3,y_1,y_2,y_3)\in \mathbf{R}^6$ to the 
system of polynomial equations
\begin{equation}\label{syst96}
d_{mn}(x_1,x_2,x_3,y_1,y_2,y_3)=0,\quad 0\le m,n\le 2,
\end{equation}
that satisfies $x_3\ge 0$ and $y_3>0$ will also satisfy conditions (i) and 
(ii) above. It turns out that if we let $x_3=0$ then system~\eqref{syst96} 
may be reduced to a system of five linear equations in the variables 
$x_1,x_2,y_1,y_2,y_3$ which admits a unique solution. We arrive in this way 
at the following solution to system~\eqref{syst96}:
\begin{equation}\label{solsyst96}
\begin{split}
&x_1=\alpha_1,\quad x_2=\alpha_2,\quad x_3=0,\quad 
y_1=\frac{3\sqrt{3}}{2r(2-3r)},\\
&y_2=-\frac{3\sqrt{3}(12r^2+8r+3)}{8r^2(3r-2)^2},
\quad y_3=C:=\frac{3}{4r(2-3r)},
\end{split}
\end{equation}
where $r$, $\alpha_1$ and $\alpha_2$ are as in~\eqref{var1}. Henceforth we 
assume that the values of the parameters $x_1,x_2,x_3,y_1,y_2,y_3$ are those 
listed in~\eqref{solsyst96}. In particular, this implies that 
$P_a(z;x_1,x_2,x_3)$ is the same as the polynomial $p_a(z)$ given 
in~\eqref{var2}. Note that by~\eqref{var1} one has 
$\Delta(p_a,z^4+z)=\mathcal{O}(a)$ for all small positive $a$. To simplify 
the notations, the zeros of the polynomial $Q_a(z;y_1,y_2,y_3)$ will be 
denoted by $\omega_1(a)$, $\omega_2(a)$ and $\overline{\omega_2(a)}$, 
respectively. Then~\eqref{omegas} and~\eqref{solsyst96} imply that
\begin{equation}\label{concl}
|\omega_1(a)-a|=|\omega_2(a)-a|=r+Ca^2.
\end{equation}

It is now practically clear that the desired conclusion should follow from 
the Newton-Raphson algorithm. We check this in a rigorous way by using 
the following simple observation.

\begin{lemma}\label{NR}
Let $R$ be a complex polynomial of degree $d\ge 2$. If $w\in \mathbf{C}$ is 
such that $R'(w)\neq 0$ then there exists a zero $z$ of $R$ such that
$$|w-z|\le d\left|\frac{R(w)}{R'(w)}\right|.$$
\end{lemma}

\begin{proof}
If $R(w)=0$ there is nothing to prove. Otherwise, we let 
$z_1,\ldots,z_d$ denote the zeros of $R$. Then
$$\sum_{i=1}^{d}\frac{1}{|w-z_i|}\ge \left|\sum_{i=1}^{d}\frac{1}{w-z_i}\right|
=\left|\frac{R'(w)}{R(w)}\right|,$$
which proves the lemma.
\end{proof}

Let $w_1(a)$, $w_2(a)$ and $w_3(a)$ denote the critical points of the 
polynomial $p_a(z)$. Since $p_a(z)\rightarrow z^4+z$ as $a\rightarrow 0$, 
we may label these critical points so that $w_1(a)\in \mathbf{R}$, 
$\Im (w_2(a))>0$ and $w_3(a)=\overline{w_2(a)}$ if $a$ is positive and 
sufficiently small. A straightforward computation shows that
\begin{equation*}
p_{a}''(\omega_1(a))=\frac{3}{r}+\mathcal{O}(a)\,\text{ and }\,
p_{a}''(\omega_2(a))=\overline{p_{a}''\left(\overline{\omega_2(a)}\right)}
=\frac{3e^{\frac{2\pi i}{3}}}{r}+\mathcal{O}(a),
\end{equation*}
which combined with Lemma~\ref{NR} and condition (ii) implies that
$$w_1(a)=\omega_1(a)+\mathcal{O}(a^3)\,\text{ and }\,w_2(a)=\overline{w_3(a)}
=\omega_2(a)+\mathcal{O}(a^3).$$
The theorem is now a consequence of~\eqref{concl}. 
\end{proof}

As one may expect, the second order variational method that we used in the 
proof of Theorem~\ref{secord4} works even in more general cases. However, for 
higher degrees the procedure described above requires an increasingly large 
amount of computations. Tedious as they may be when done by hand, for small 
degrees these computations become relatively easy if one uses for instance a 
Maple computer program. Indeed, such a program has considerably simplified 
our task in the course of proving Theorem~\ref{secord4} and it also allowed 
us to obtain a similar result for the polynomial $z^5+z$. In order to formulate
this result we need to introduce some additional notations: let $a\in [0,1]$
and set 
\begin{equation*}
\begin{split}
&s=\left|z^5+z\right|_0=\sqrt[4]{\frac{1}{5}},\quad
\beta=\frac{2\sqrt{2}s^2}{1-2s^2},\quad
\gamma=\frac{4\sqrt{2}}{5s(1-2s^2)},\\
&\delta=\frac{60s^4-19}{50s^2(2s^2-1)^2},\quad 
K=\frac{2}{5s(1-2s^2)}\approx 5.665658792,\\
&\eta(a)=\exp{\left[i\!\left(\frac{\pi}{4}+\beta a\right)\right]},\quad 
\chi_1(a)=a+\left(s+Ka^2\right)
\exp{\left[i\!\left(\frac{\pi}{4}+\gamma a+\delta a^2\right)\right]},\\
&\chi_2(a)=a+\left(s+Ka^2\right)
\exp{\left[i\!\left(\frac{3\pi}{4}+\gamma a-\delta a^2\right)\right]}.
\end{split}
\end{equation*}
We use these quantities in order to define the following one-parameter 
families of polynomials:
\begin{eqnarray*}
&&q_a(z):=(z-a)(z-\eta(a))(z-i\eta(a))\!\left[z-\overline{\eta(a)}\right]
\!\left[z+i\overline{\eta(a)}\right]\in S_5,\\
&&s_a(z):=5(z-\chi_1(a))(z-\chi_2(a))\!\left[z-\overline{\chi_1(a)}\right]
\!\left[z-\overline{\chi_2(a)}\right].
\end{eqnarray*}
Our next result is an analogue of Theorem~\ref{secord4} for the polynomial 
$z^5+z$. Its proof is similar to that of Theorem~\ref{secord4} and is 
therefore omitted.

\begin{theorem}\label{secord5}
For all sufficiently small $a>0$ one has
\begin{equation*}
\Delta(q_a,z^5+z)=\mathcal{O}(a)\,\text{ and }\,
\Delta\!\left(\frac{q'_a}{5},\frac{s_a}{5}\right)=\mathcal{O}(a^3),
\end{equation*}
so that $d(q_a)=|q_a|_a=s+Ka^2+\mathcal{O}(a^3)$. In particular, 
the polynomial $z^5+z$ is not locally maximal for Sendov's conjecture.
\hfill $\Box$
\end{theorem}

\begin{remark}\label{gencase}
In view of Lemma~\ref{countlm3} and Theorems~\ref{secord4}--\ref{secord5} it 
seems reasonable to conjecture that polynomials of 
the form $z^n+z$ or, more generally, the 0-maximal polynomials 
of degree $n\ge 2$ given in Theorem~\ref{class0max} are not local maxima for 
the function $d$. One could for instance try to find an algorithmic proof of 
this conjecture based on the 
second order variational method that we just described. However, we believe 
that the results of this section are enough to emphasize both the limitations 
of what can be achieved through first order variational methods and the 
necessity of using
higher order approximations when dealing with problems such as Sendov's 
conjecture. In particular, these results reinforce our earlier statement 
that methods such as those used in \cite{S1} and \cite{S2} can only lead to 
already known partial results and are unlikely to be 
successful in the general case.
\end{remark}

As we already pointed out in the introduction, the fact that the function 
$d\circ \tau$ fails to be (logarithmically) plurisubharmonic in the polydisk 
$\bar{D}^n$ accounts 
for many of the difficulties in studying locally maximal polynomials for 
Sendov's conjecture. Nevertheless, determining all such polynomials is a 
natural and interesting question. The polynomials in Conjecture~\ref{send2} 
(i.~e., the polynomials which are conjectured to be extremal for Sendov's 
conjecture) are the only examples of local maxima for the function $d$ known 
so far (cf.~\cite{M3} and \cite{VZ}). We propose the following  
stronger version of Conjecture~\ref{send2}: 

\begin{conjecture}\label{send3}
Let $p\in S_n$. Then $p$ is locally maximal for Sendov's conjecture if and 
only if $p(z)=z^n+e^{i\theta}$ for some $\theta\in \mathbf{R}$.
\end{conjecture}

To the best of our 
knowledge, Conjecture~\ref{send3} has not been stated explicitly in the 
literature. As one can see from the discussion in sections 1 and 2, any 
approach based on local variational methods that were to confirm 
Conjecture~\ref{send2} would also have to confirm Conjecture~\ref{send3}.

\section{Differentiators of normal operators, majorization, and the 
geometry of polynomials}

The results known so far on Sendov's conjecture and related questions were 
almost exclusively 
obtained by analytical arguments. In this section we propose an operator 
theoretical 
approach to Conjectures~\ref{send1} and~\ref{send2} and show that these may 
be viewed as part of the more general problem of describing the 
relationships between the spectra of normal matrices and the spectra of their 
principal submatrices. We also give a geometrical characterization of the 
critical points of complex polynomials by means of multivariate majorization 
relations.

\subsection{Sendov's conjecture and spectral variations of normal operators 
and their compressions}

Let $\mathcal{H}$ be a $n$-dimensional complex Hilbert space with (unitarily 
invariant) scalar 
product $\langle\cdot,\!\cdot\rangle$ and identity operator 
$I_{\mathcal{H}}\in L(\mathcal{H})$, where $L(\mathcal{H})$ is the
set of all linear operators on $\mathcal{H}$. We may view the spectrum
$\text{Eig}(A)$ of an operator $A\in L(\mathcal{H})$ as the multiset whose 
elements are the eigenvalues of $A$, each eigenvalue occurring as many times 
as its algebraic multiplicity. The spectral radius of $A$ is then given 
by $\rho(A)=\max_{z\in \text{Eig}(A)}|z|$. Let $\mathcal{K}$ be a 
subspace of $\mathcal{H}$, $A\in L(\mathcal{H})$ and $B\in L(\mathcal{K})$. 
Following \cite{He} (see also \cite{Sc}), we define the {\em spectral 
variation} $s(A,B)$ of $A$ and $B$ to be the directed Hausdorff distance from 
$\text{Eig}(A)$ to $\text{Eig}(B)$, that is,
$$s(A,B)=\max_{z\in \text{Eig}(A)}\min_{w\in \text{Eig}(B)}|z-w|.$$

It was shown in \cite[Theorem 1.11]{P} that any $A\in L(\mathcal{H})$ has a 
so-called {\em trace vector}, i.~e., a vector $v$ that satisfies 
$v^{*}A^{k}v=\frac{1}{n}\text{tr}(A^k)$ for any non-negative integer $k$. 
It is not difficult to see that if $A\in L(\mathcal{H})$ is normal (i.~e., 
$AA^{*}=A^{*}A$) and 
$(\mathbf{e}_1,\ldots,\mathbf{e}_n)$ is an orthonormal basis of 
$\mathcal{H}$ consisting of 
eigenvectors for $A$ then $\mathbf{v}_{n}:=\frac{1}{\sqrt{n}}(\mathbf{e}_1+
\ldots+\mathbf{e}_n)$ is a 
trace vector of $A$. Let $P$ denote the orthoprojection on the subspace 
$\mathcal{K}:=\mathbf{v}_{n}^{\perp}$ of $\mathcal{H}$ and define the $P$-{\em 
compression} of $A$ to be 
the operator $A'=PAP|_{\mathcal{K}}\in L(\mathcal{K})$ (cf.~\cite{D}). Then 
one can show the following result.

\begin{proposition}\label{normal}
Let $A\in L(\mathcal{H})$ be a normal operator with spectrum 
$\text{{\em Eig}}(A)=\{z_1,\ldots,z_n\}$ and let $A'$ be as above. Then
$$\det(A'-zI_{\mathcal{K}})=\frac{1}{n}\!
\left[\sum_{i=1}^{n}\frac{1}{z_i-z}\right]\det
(A-zI_{\mathcal{H}})=-\frac{1}{n}\frac{d}{dz}
\left[{\det(A-zI_{\mathcal{H}})}\right]$$
for any $z\notin \text{{\em Eig}}(A)$.
\end{proposition}
\begin{proof}
Let $(\mathbf{v}_1,\ldots,\mathbf{v}_{n-1})$ be an orthonormal basis of 
$\mathcal{K}$. Then the 
matrix representation of $A'$ in the basis $(\mathbf{v}_1,\ldots,
\mathbf{v}_{n-1})$ is given by 
the $(n-1)\times (n-1)$ upper left-hand principal submatrix of the matrix 
representation of $A$ in the orthonormal basis $(\mathbf{v}_1,\ldots,
\mathbf{v}_{n-1},\mathbf{v}_n)$ of 
$\mathcal{H}$. For any $z\notin \text{{Eig}}(A)$ the $(n,n)$ entry of the 
matrix representation of the (normal) operator $(A-zI_{\mathcal{H}})^{-1}$ in 
the basis $(\mathbf{v}_1,\ldots,\mathbf{v}_{n-1},\mathbf{v}_n)$ is given on 
the one hand by
\begin{eqnarray*}
\left\langle (A-zI_{\mathcal{H}})^{-1}\mathbf{v}_n,\mathbf{v}_n\right\rangle 
&=&
\left\langle \sum_{i=1}^{n}(A-zI_{\mathcal{H}})^{-1}\langle \mathbf{v}_n,
\mathbf{e}_i\rangle \mathbf{e}_i,
\sum_{i=1}^{n}\langle \mathbf{v}_n,\mathbf{e}_i\rangle \mathbf{e}_i
\right\rangle \\
&=&\left\langle \sum_{i=1}^{n}\frac{\langle \mathbf{v}_n,\mathbf{e}_i\rangle}
{z_i-z}\mathbf{e}_i,
\sum_{i=1}^{n}\langle \mathbf{v}_n,\mathbf{e}_i\rangle \mathbf{e}_i
\right\rangle=\sum_{i=1}^{n}
\frac{|\langle \mathbf{v}_n,\mathbf{e}_i\rangle|^{2}}{z_i-z}\\
&=&\frac{1}{n}
\sum_{i=1}^{n}\frac{1}{z_i-z}.
\end{eqnarray*}
On the other hand, by Cramer's rule the $(n,n)$ entry of the matrix 
representation of $(A-zI_{\mathcal{H}})^{-1}$ in 
the basis $(\mathbf{v}_1,\ldots,\mathbf{v}_{n-1},\mathbf{v}_n)$ of 
$\mathcal{H}$ is given by the cofactor of the $(n,n)$ 
entry of the matrix representation of $(A-zI_{\mathcal{H}})$ in the same 
basis. Thus 
$$\left\langle (A-zI_{\mathcal{H}})^{-1}\mathbf{v}_n,\mathbf{v}_n\right\rangle=
\frac{\det(A'-zI_{\mathcal{K}})}{\det(A-zI_{\mathcal{H}})},$$
which proves the proposition.
\end{proof}

\begin{remark}\label{normcompr}
The $P$-compression $A'$ of the normal operator $A$ is not necessarily normal 
itself. As it was shown in~\cite[Proposition 3.4]{P}, this occurs if and only 
if the eigenvalues of $A$ are collinear complex numbers.
\end{remark}  

By Proposition~\ref{normal}, up to a factor $-\frac{1}{n}$ the characteristic 
polynomial of the $P$-compression $A'$ of $A$ coincides with the derivative of 
the characteristic polynomial of $A$. For this reason, the orthoprojection 
$P$ is called a {\em differentiator} of $A$ (cf.~\cite{D} 
and~\cite[Definition 1.3]{P}). Clearly, for any complex polynomial $p(z)$ of 
degree $n\ge 2$ 
there exists a normal operator $A\in L(\mathcal{H})$ such that the 
characteristic polynomial of $A$ is $(-1)^{n}p(z)$. Indeed, 
if $p(z)=\prod_{i=1}^{n}(z-z_i)$ then one may just consider an
operator whose matrix representation in a given orthonormal basis 
of $\mathcal{H}$ is $\text{diag}(z_1,\ldots,z_n)$. We deduce that 
Conjectures 1 and 2 can 
be formulated exclusively in terms of the spectra of a normal operator and 
its compression to the orthogonal complement of a 
trace vector:

\begin{conjecture}\label{operator}
Let $A$ be a normal operator on a complex Hilbert space of dimension $n\ge 2$ 
and let $A'$ denote the compression of $A$ to the orthogonal complement of a 
trace vector of $A$. Then $s(A,A')\le \rho(A)$
and equality occurs if and only if either $A=0$ or $A$ has simple eigenvalues 
and $A^n$ is a non-zero scalar operator.
\end{conjecture}

Note that by the results of \cite{BX} we know that 
Conjecture~\ref{operator} is true for normal operators on a complex Hilbert 
space of dimension $n\le 8$. A natural question that arises in this context 
is whether there exists an analogue of Conjecture~\ref{operator} for the 
directed Hausdorff distance from $\text{Eig}(A')$ to $\text{Eig}(A)$. The 
following theorem gives a complete answer to the latter question. 

\begin{theorem}\label{converse}
Let $A$ be a normal operator on a complex Hilbert space of dimension $n\ge 2$ 
and let $A'$ denote the compression of $A$ to the orthogonal complement of a 
trace vector of $A$. Then $s(A',A)\le \rho(A)$
and equality occurs if and only if $\text{{\em tr}}(A)=0$ and 
$\min_{z\in \text{{\em Eig}}(A)}|z|=\rho(A)$.
\end{theorem}

Theorem~\ref{converse} is an immediate consequence of the following 
more general result:

\begin{lemma}\label{simplelemma}
Let $M$ and $N$ be finite point sets in the complex plane such that $M$ is 
contained in the closed convex hull of the points in $N$. Then 
$$\max_{w\in M}\min_{z\in N}|w-z|\le \max_{z\in N}|z|$$ 
and equality occurs if and only if $0\in M$ and 
$\min_{z\in N}|z|=\max_{z\in N}|z|$.
\end{lemma}
\begin{proof}
It is geometrically clear that $\min_{z\in N}|w-z|^2\le
\max_{z\in N}|z|^2-|w|^2$ for any $w\in M$. This implies that if 
$0\neq w\in M$ then $\min_{z\in N}|w-z|<\max_{z\in N}|z|$.
\end{proof}

A characterization of trace vectors and 
differentiators of arbitrary operators on finite-dimensional Hilbert spaces 
was given in~\cite[Theorem 2.5]{P}. Using this characterization we can 
strengthen Proposition~\ref{normal} in the following way:

\begin{proposition}\label{ONtrace}
Let $A\in L(\mathcal{H})$ be a normal operator with spectrum 
$\text{{\em Eig}}(A)=\{z_1,\ldots,z_n\}$ and let 
$(\mathbf{e}_1,\ldots,\mathbf{e}_n)$ be an orthonormal basis of $\mathcal{H}$ 
consisting of eigenvectors for $A$. Set
$$\mathbf{v}_i=\frac{1}{\sqrt{n}}\sum_{j=1}^{n}\eta^{ij}\mathbf{e}_j,\quad 
1\le i\le n,$$
where $\eta=\exp\!\left(\frac{2\pi \sqrt{-1}}{n}\right)$. Then 
$(\mathbf{v}_1,\ldots,\mathbf{v}_n)$ is an orthonormal basis of $\mathcal{H}$ 
consisting of trace vectors for $A$. Thus, if $P_i$ denotes the 
orthoprojection on the subspace 
$\mathcal{K}_i:=\mathbf{v}_i^{\perp}=\mathbf{C}\text{-span}
\{\mathbf{v}_1,
\ldots,\hat{\mathbf{v}}_i,\ldots,\mathbf{v}_n\}$ of $\mathcal{H}$ and 
$A_i=P_iAP_i|_{\mathcal{K}_i}\in L(\mathcal{K}_i)$ is the $P_i$-compression of 
$A$ then
$$\det(A_i-zI_{\mathcal{K}_i})=\frac{1}{n}\!
\left[\sum_{i=1}^{n}\frac{1}{z_i-z}\right]\det
(A-zI_{\mathcal{H}})=-\frac{1}{n}\frac{d}{dz}
\left[{\det(A-zI_{\mathcal{H}})}\right]$$
for any $z\notin \text{{\em Eig}}(A)$.
\end{proposition}
\begin{proof}
It is easy to check that $(\mathbf{v}_1,\ldots,\mathbf{v}_n)$ is an 
orthonormal basis of $\mathcal{H}$. The same computations as in the proof 
of Proposition~\ref{normal} show that one the one hand
\begin{equation*}
\left\langle (A-zI_{\mathcal{H}})^{-1}\mathbf{v}_i,\mathbf{v}_i\right\rangle 
=\sum_{j=1}^{n}
\frac{|\langle \mathbf{v}_i,\mathbf{e}_j\rangle|^{2}}{z_j-z}\\
=\frac{1}{n}
\sum_{j=1}^{n}\frac{1}{z_j-z}
\end{equation*}
and on the other hand
$$\left\langle (A-zI_{\mathcal{H}})^{-1}\mathbf{v}_i,\mathbf{v}_i\right\rangle=
\frac{\det(A_i-zI_{\mathcal{K}_i})}{\det(A-zI_{\mathcal{H}})}$$
for $1\le i\le n$. By~\cite[Theorem 2.5]{P} this is the same as saying that 
the vectors $\mathbf{v}_i$, $1\le i\le n$, are trace vectors of the operator 
$A$.
\end{proof}

The following consequence of Proposition~\ref{ONtrace} should prove useful 
for studying the geometry of the zeros and critical points of complex
polynomials by operator theoretical methods.

\begin{corollary}\label{matrpoly}
Let $p$ be a monic complex polynomial of degree $n\ge 2$. Then there exists a 
normal $n\times n$ complex matrix $\mathbf{A}$ with the following properties:
$$\det(\mathbf{A}-z\mathbf{I}_n)=(-1)^np(z)\text{ and }
\det(\mathbf{A}_{[i]}-z\mathbf{I}_{n-1})=\frac{(-1)^{n-1}}{n}p'(z),
\,1\le i\le n,$$
where $\mathbf{A}_{[i]}$ denotes the degeneracy one principal submatrix of 
$\mathbf{A}$ which is obtained by deleting the $i$-th row 
and $i$-th column from $\mathbf{A}$. 
\end{corollary}
\begin{proof}
Let $p(z)=\prod_{i=1}^{n}(z-z_i)$ and consider the
operator $A$ whose matrix representation in a given orthonormal basis 
$(\mathbf{e}_1,\ldots,\mathbf{e}_n)$ of 
$\mathcal{H}$ is $\mathbf{B}=\text{diag}(z_1,\ldots,z_n)$. Define the 
following unitary $n\times n$ matrix:
$$\mathbf{U}=(u_{ij}),\text{ where } u_{ij}=\frac{\eta^{ij}}{\sqrt{n}},\,
1\le i,j\le n,\text{ and }\eta=\exp\!\left(\frac{2\pi \sqrt{-1}}{n}\right),$$
and set $\mathbf{A}=\mathbf{U}^{*}\mathbf{B}\mathbf{U}$. Then $\mathbf{A}$ is 
the matrix representation of the operator $A$ in the orthonormal basis 
$(\mathbf{v}_1,\ldots,\mathbf{v}_n)$ of $\mathcal{H}$ which was constructed 
in Proposition~\ref{ONtrace}. Moreover, 
$\mathbf{A}_{[i]}$ is precisely the matrix representation of the 
$P_i$-compression $A_i$ of $A$ in the orthonormal basis $(\mathbf{v}_1,
\ldots,\hat{\mathbf{v}}_i,\ldots,\mathbf{v}_n)$ of $\mathbf{v}_i^{\perp}$ for 
$1\le i\le n$. It follows from Proposition~\ref{ONtrace} that 
$\mathbf{A}$ must have the desired properties.
\end{proof}

The results that we have just presented make use of the ideas 
developed by Davis in \cite{D} and Pereira in \cite{P}. In conjunction 
with methods of majorization theory, the operator theoretical tools 
introduced in \cite{P} were 
the key to Pereira's main results, namely some beautiful solutions to the 1947 
conjecture of de Bruijn and Springer and the 1986 conjecture of Schoenberg. 
Similar ideas were used by Malamud in \cite{Ma}, where he not only 
proved these same two conjectures -- independently and almost at the same 
time as Pereira -- but he also obtained a remarkable generalization of the de 
Bruijn-Springer conjecture (see section 3.2 below). The methods used in 
\cite{Ma} and \cite{P} seem to be particularly well suited for studying 
extremal problems for which the loci of the zeros of extremal polynomials are 
(conjectured to be) lines in the complex plane (cf.~Remark~\ref{normcompr}). 
Notwithstanding, the ideas of 
\cite{Ma} and \cite{P} are a new and powerful approach to the study of the 
geometry of polynomials. We believe that the results and the setting 
developed in this section should prove useful for 
investigating Conjecture~\ref{operator}, that is, Sendov's conjecture 
formulated in terms of normal operators.

\subsection{Majorization and a characterization of the critical points of 
complex polynomials}

Although useful in many contexts, the geometrical information contained in the 
Gauss-Lucas theorem is hardly sufficient for dealing with problems such as 
Sendov's conjecture. This is mainly because of the implicit nature of the 
relations between the zeros and the critical points, which actually accounts 
for many of the difficulties in studying the geometry of polynomials. 
Describing these relations geometrically and as explicitly as 
possible would be helpful for a great many questions in 
the analytic theory of polynomials. Indeed, most of these questions relate in 
various ways to the following fundamental problem:

\begin{problem}\label{fund}
Let $z_1,\ldots,z_n$ and $w_1,\ldots,w_{n-1}$ be $2n-1$ points in the complex 
plane, where $n\ge 2$. Find necessary and sufficient geometric conditions 
in order for $w_1,\ldots,w_{n-1}$ to be the critical 
points of the polynomial $\prod_{i=1}^{n}(z-z_i)$.
\end{problem} 

Naturally, there may be several possible ways of answering Problem~\ref{fund}
depending on the context in which one places this problem and the approach 
that one uses. In this section we build on the results of \cite{Ma} in 
order to give an answer to  
Problem~\ref{fund} in terms of majorization, which is a fundamental concept 
in the theory of inequalities as well as matrix analysis and operator theory 
(see \cite{MO}).

Let us start with the following definition.

\begin{definition}\label{stoc}
An $m\times n$ matrix $R=(r_{ij})$ with $m\le n$ is called {\em rectangularly 
stochastic} if it satisfies the following conditions:
$$r_{ij}\ge 0,\,\sum_{j=1}^{n}r_{ij}=1\text{ and }\sum_{i=1}^{m}r_{ij}
=\frac{m}{n}\text{ for }1\le i\le m,\,
1\le j\le n.$$
\end{definition}

Note that if $m=n$ then Definition~\ref{stoc} reduces to the usual definition 
of a doubly stochastic (or bistochastic) matrix. In what follows we shall 
identify an $m$-tuple $X=(x_1,\ldots,x_m)^t$ consisting of vectors in 
$\mathbf{R}^k$ with an $m\times k$ real matrix in the obvious manner.

The following theorem is 
a fundamental result in the theory of multivariate majorization
which is originally due to Sherman (\cite{Sher}).

\begin{theorem}\label{sher}
Let $X=(x_1,\ldots,x_m)^t$ and $Y=(y_1,\ldots,y_n)^t$ be an unordered 
$m$-tuple and $n$-tuple, respectively, of vectors in 
$\mathbf{R}^k$, where $n\ge m$ and $k\ge 1$. Then the following conditions 
are equivalent:
\begin{itemize}
\item[(i)] For any convex function $f:\mathbf{R}^k\rightarrow \mathbf{R}$ one 
has $\dfrac{\sum_{i=1}^{m}f(x_i)}{m}\le \dfrac{\sum_{i=1}^{n}f(y_i)}{n}$.
\item[(ii)] There exists a rectangularly stochastic $m\times n$ matrix $R$ 
such that $\tilde{X}=R\tilde{Y}$, where $\tilde{X}$ is an $m\times k$ matrix 
and $\tilde{Y}$ is an $n\times k$ matrix obtained by some (and then any) 
ordering of the vectors in $X$ and $Y$.
\end{itemize}
If (i) and (ii) are satisfied then we say that $X$ is {\em majorized} by 
$Y$ and write $X\prec Y$.
\end{theorem} 

If $X$ and $Y$ consist of complex numbers then one defines the 
majorization relation $X\prec Y$, when 
appropriate, by identifying 
$\mathbf{C}$ with $\mathbf{R}^2$ in the above theorem.

\begin{remark}
If $m=n$ then Birkhoff's theorem (\cite[Theorem A.2]{MO}) implies 
that the majorization relation $\prec$ defines a partial ordering on the set 
of unordered $n$-tuples of vectors in $\mathbf{R}^k$. The case $k=1$ and $m=n$ 
is originally due to Schur and to Hardy, Littlewood and P\'olya and is often 
referred to as classical majorization. We refer to \cite{MO} for further 
details, including a discussion of the geometrical meaning of classical 
majorization.
\end{remark}

\begin{remark}
Surprisingly, Theorem~\ref{sher} was long assumed to be an open problem and 
it does not appear in \cite{MO}, which is the 
definite reference on majorization theory. This may explain why Sherman's 
result has been rediscovered several times.
\end{remark}

Any geometric description of the critical points of a complex polynomial 
is rooted in, and so it should also reflect, some algebraic characterization 
of these points. One such algebraic characterization is given by the following
lemma.

\begin{lemma}\label{alg}
Let $p(z)=\prod_{i=1}^{n}(z-z_i)$ and $q(z)=n\prod_{j=1}^{n-1}(z-w_j)$. 
The following conditions are equivalent:
\begin{itemize}
\item[(i)] $q(z)=p'(z)$ for any $z\in \mathbf{C}$.
\item[(ii)] For any $\alpha\in \mathbf{C}$ and $1\le k\le n-1$
one has
\begin{equation*}
\frac{1}{\binom{n-1}{k}}\sum_{1\le i_1<\ldots<i_k\le n-1}\prod_{j=1}^{k}
(w_{i_j}-\alpha)=\frac{1}{\binom{n}{k}}\sum_{1\le i_1<\ldots<i_k\le n}
\prod_{j=1}^{k}(z_{i_j}-\alpha).
\end{equation*}
\item[(iii)] The relations in {\em (ii)}$\!$ are valid for $\alpha=0$.
\end{itemize}
\end{lemma}

\begin{proof}
Let $\alpha\in \mathbf{C}$. Then (i) is obviously equivalent to 
$q^{(n-1-k)}(\alpha)=p^{(n-k)}(\alpha)$ for $1\le k\le n-1$. The
lemma follows from elementary computations. 
\end{proof}

By combining \cite[Theorem 4.6]{Ma} with Theorem~\ref{sher} and 
Lemma~\ref{alg} we get the following  
answer to Problem~\ref{fund} in terms of
multivariate majorization relations and convex functions: 

\begin{theorem}\label{main}
Let $p(z)=\prod_{i=1}^{n}(z-z_i)$, $q(z)=n\prod_{j=1}^{n-1}(z-w_j)$ and
 $\alpha\in \mathbf{C}$. For $1\le k\le n-1$
define the following unordered $\tbinom{n}{k}$-tuple and 
$\tbinom{n-1}{k}$-tuple, respectively, of complex numbers:
\begin{eqnarray*}
Z(\alpha,k)&=&\bigg(\prod_{j=1}^{k}(z_{i_j}-\alpha)\bigg)^{t}_{1\le i_1<\ldots
<i_k\le n},\\
W(\alpha,k)&=&\bigg(\prod_{j=1}^{k}(w_{i_j}-\alpha)\bigg)^{t}_{1\le 
i_1<\ldots<i_k\le n-1}.
\end{eqnarray*}
Then the following conditions are equivalent:
\begin{itemize}
\item[(i)] $q(z)=p'(z)$ for any $z\in \mathbf{C}$.
\item[(ii)] For any $\alpha\in \mathbf{C}$ and $1\le k\le n-1$ one 
has $W(\alpha,k)\prec Z(\alpha,k)$.
\item[(iii)] For $1\le k\le n-1$ one 
has $W(0,k)\prec Z(0,k)$.
\item[(iv)] For any convex function $f:\mathbf{C}\rightarrow \mathbf{R}$, 
$\alpha \in \mathbf{C}$ and $1\le k\le n-1$ one has
\begin{equation*}
\frac{1}{\binom{n-1}{k}}\sum_{1\le i_1<\ldots<i_k\le n-1}f\!
\bigg(\!\prod_{j=1}^{k}(w_{i_j}-\alpha)\!\bigg)\le \frac{1}{\binom{n}{k}}
\sum_{1\le i_1<\ldots<i_k\le n}f\!\bigg(\!\prod_{j=1}^{k}(z_{i_j}-\alpha)
\!\bigg).
\end{equation*}
\item[(v)] The relations in {\em (iv)}$\!$ are valid for $\alpha=0$.
\end{itemize}
\end{theorem} 

The implication (i) $\Rightarrow$ (iv) with $k=1$ and $\alpha=0$ was
originally conjectured by de Bruijn and Springer in \cite{BS} and is in 
itself a powerful 
generalization of the Gauss-Lucas theorem (see also \cite[Theorem 5.4]{P} for 
a proof of the de Bruijn-Springer conjecture). Unlike classical majorization, 
the geometrical meaning of multivariate majorization -- in our case, 
majorization in $\mathbf{R}^2$ -- is not yet fully understood. In view of the 
remarkable achievements of \cite{Ma} and \cite{P}, further clarifying the 
geometrical information encoded in Theorem~\ref{main} should prove very useful
for studying the geometry of polynomials in general and Problem~\ref{fund} 
and Sendov's conjecture in particular.

\subsection{The spectra of principal submatrices of normal matrices}

It is interesting to note that the method used in the proof of 
Proposition~\ref{normal} leads actually to an analogue of the Gauss-Lucas 
theorem for normal matrices. Indeed, 
let $n\ge 2$ and let $A$ be a normal $n\times n$ complex matrix with 
spectrum $\text{Eig}(A)=\{z_1,\ldots,z_n\}$. As in section 3.1, for 
$1\le i\le n$ we denote 
by $A_{[i]}$ the degeneracy one principal submatrix of $A$ obtained by 
deleting the $i$-th row and $i$-th column from $A$. Let $I_n$ be the 
$n\times n$ identity matrix and 
$(\textbf{u}_1,\ldots,\textbf{u}_n)$ 
be the standard basis of $\mathbf{C}^n$ viewed as a complex Hilbert 
space with standard scalar product. Now fix an orthonormal basis 
$(\mathbf{e}_1,\ldots,\mathbf{e}_n)$ of $\mathbf{C}^n$ consisting of 
eigenvectors for $A$. Using the same arguments as in the proof of 
Proposition~\ref{normal} we get the following theorem of Gauss-Lucas type for 
normal matrices:

\begin{theorem}\label{GL}
With the above notations one has
$$\frac{\det(A_{[i]}-zI_{n-1})}{\det(A-zI_n)}
=\sum_{j=1}^{n}\frac{|\langle \mathbf{u}_i,\mathbf{e}_j \rangle|^2}{z_j-z}$$
for $z\in \mathbf{C}\setminus \text{{\em Eig}}(A)$ and $1\le i\le n$.
In particular, $\text{{\em Eig}}(A_{[i]})$ is contained in the convex hull 
of $\text{{\em Eig}}(A)$ for any $i\in \{1,2,\ldots,n\}$.\hfill $\Box$
\end{theorem}

Let $\det(A-zI_n)=\prod_{k=1}^{m}(z_k-z)^{n_k}$, where $\sum_{k=1}^{m}n_k=n$ 
and the $z_k$'s denote the distinct eigenvalues of $A$ with 
multiplicities $n_k\ge 1$, $1\le k\le m$, respectively. Then the $z_k$'s 
are eigenvalues of each of the matrices $A_{[i]}$, $1\le i\le n$, with 
multiplicities at least $n_k-1$, $1\le k\le m$, respectively. Thus,
 there exist complex numbers $w_{j}^{(i)}$, 
$1\le i\le n$, $1\le j\le m-1$, such that
$$\det(A_{[i]}-zI_{n-1})=\prod_{j=1}^{m-1}
(w_{j}^{(i)}-z)\prod_{k=1}^{m}(z_k-z)^{n_k-1}.$$
By Theorem~\ref{GL} we get
$$\frac{\prod_{j=1}^{m-1}(w_{j}^{(i)}-z_k)}
{\prod_{j=1,\,j\neq k}^{m}(z_j-z_k)}=n_k|\langle \mathbf{u}_i,\mathbf{e}_k 
\rangle|^2,\quad 1\le k\le m,\,1\le i\le n,$$
which proves the following proposition:

\begin{proposition}\label{int}
With the above notations one has
\begin{equation}\label{inter}
\frac{\prod_{j=1}^{m-1}(w_{j}^{(i)}-z_k)}
{\prod_{j=1,\,j\neq k}^{m}(z_j-z_k)}\ge 0
\end{equation}
for all $1\le k\le m$ and $1\le i\le n$.\hfill $\Box$
\end{proposition}

Note that if $A$ is Hermitian then \eqref{inter} simply means that 
for each $i\in \{1,\ldots,n\}$ the sequences of real numbers 
$w_{1}^{(i)},\ldots,w_{m-1}^{(i)}$ and $z_1,\ldots,z_m$ are interlacing, which 
is the well-known Cauchy-Poincar\'e separation theorem for self-adjoint 
matrices (cf., e.~g., \cite{D}, where one may also find a converse to the 
Cauchy-Poincar\'e theorem due to Ky Fan and Pall). 
Proposition~\ref{int} may therefore be viewed as an analogue of the 
Cauchy-Poincar\'e theorem for normal matrices.

Proposition~\ref{int} is a slightly more general version of the necessity 
part of Proposition 3.1 in \cite{Ma}, where a partial converse to this 
result was obtained. Note though that the conditions stated in 
{\em loc.~cit.~}are 
actually valid only in the generic case when $A$ has simple eigenvalues. 
Indeed, the arguments used above show that Proposition 3.1 in \cite{Ma} 
should be stated as follows in the general case: 

\begin{proposition}
Let $\lambda_1,\ldots,\lambda_n$ and $\mu_1,\ldots,\mu_{n-1}$ be 
$2n-1$ complex numbers. The conditions
$$\lim_{\lambda\rightarrow \lambda_k}\!\left[\frac{\prod_{j=1}^{n-1}
(\mu_j-\lambda)}{\prod_{j=1,\,j\neq k}^{n}(\lambda_j-\lambda)}\right]\ge 0,
\quad 1\le k\le n,$$
are valid if and only if there exists a (not necessarily unique) normal 
$n\times n$ matrix $A$ such that $\text{{\em Eig}}(A)=
\{\lambda_1,\ldots,\lambda_n\}$
and $\text{{\em Eig}}(A_{[n]})=\{\mu_1,\ldots,\mu_{n-1}\}$.
\end{proposition}

Obviously, the critical points of the characteristic polynomial of a normal 
matrix satisfy both the ``Gauss-Lucas property'' (Theorem~\ref{GL}) and 
-- by Lemma~\ref{alg} -- the ``interlacing property'' (Proposition~\ref{int}).
In light of Theorem~\ref{GL} and Proposition~\ref{int}, we may view
the spectra of the deficiency one principal submatrices of a normal 
matrix $A$ as generalizations of the set of critical points of the 
characteristic polynomial of $A$. As already noted in section 3.1, any monic 
complex polynomial 
$p$ of degree $n$ may be realized -- up to a factor $(-1)^n$ -- as the 
characteristic polynomial of a normal $n\times n$ matrix $A$. Since $A$
has a differentiator, the polynomial $p'$ may itself be
realized -- up to a factor $\frac{(-1)^{n-1}}{n}$ and possibly also 
a unitary change of orthonormal basis -- as the characteristic polynomial of a 
(not necessarily normal) deficiency one principal 
submatrix of $A$. We conclude that
Problem~\ref{fund} and the study of the 
geometry of zeros and critical points of complex polynomials 
-- including Sendov's conjecture -- may be seen as part of the following more 
general problem:

\begin{problem}\label{norm}
Let $A$ be a normal $n\times n$ complex matrix. Describe the geometrical 
properties and the mutual location of the sets $\text{Eig}(A_{[i]})$ and 
$\text{Eig}(A)$ for $1\le i\le n$.
\end{problem}

The example $A=\text{diag}(1,1,\ldots,1,-1)$ and 
$A_{[n]}=\text{diag}(1,1,\ldots,1)$ shows that the directed Hausdorff distance
from
$\text{Eig}(A)$ to $\text{Eig}(A_{[n]})$ can actually be equal to the 
diameter of the convex hull of $\text{Eig}(A)$. Thus, the ``Gauss-Lucas 
theorem for normal matrices'' (Theorem~\ref{GL}) provides in fact an optimal 
upper bound 
for this directed Hausdorff distance, so that there is no apparent analogue 
of Sendov's 
conjecture for Problem~\ref{norm}. Nevertheless, in view of Theorem~\ref{main}
it seems natural to look for an answer to this problem in terms of (weak) 
majorization relations between eigenvalues (cf.~\cite{MO}) and inequalities of 
de Bruijn-Springer type. These 
questions will be addressed in a forthcoming paper.

\medskip

\noindent
{\bf Acknowledgements.} The author would like to thank Harold Shapiro for 
his interest in this work and Arne Meurman for useful suggestions.}

\end{document}